\documentclass[final,5p,times,twocolumn]{elsarticle}

\usepackage{amssymb}
\usepackage{amsmath}
\usepackage{amsfonts}
\usepackage{xcolor}
\usepackage{mathrsfs}
\usepackage{amsthm}
\usepackage{float}
\usepackage{xr}
\newcommand{\NY}[1]{{\color{black} #1}}
\newcommand{\NYY}[1]{{\color{black} #1}}
\newcommand{\NYYY}[1]{{\color{black} #1}}
\newcommand{\NYb}[1]{{\color{black} #1}}

\newcommand{\GY}[1]{{\color{black} #1}}

\begin{document}
	
	\begin{frontmatter}
		
		\title{
			Inverse problems of inhomogeneous fracture toughness using
			phase-field models
		}

		\author[1]{Yueyuan Gao\corref{cor1}}
		\ead{yueyuan.gao@es.hokudai.ac.jp}
		
		\author[2,3]{Natsuhiko Yoshinaga}
		\ead{yoshinaga@tohoku.ac.jp}

		\cortext[cor1]{Corresponding author}
		
		\affiliation[1]{organization={Laboratory of Mathematical Modeling, Research Institute for Electronic Science},
			addressline={Hokkaido University},
			postcode={060-0812},
			city={Sapporo},
			country={Japan}
		}

		\affiliation[2]{organization={WPI - Advanced Institute for Materials Research, Tohoku University},
			addressline={Katahira 2-1-1},
			postcode={980-8577},
			city={Sendai},
			country={Japan}
		}
		
		\affiliation[3]{organization={MathAM-OIL, AIST},
			addressline={Katahira 2-1-1},
			postcode={980-8577},
			city={Sendai},
			country={Japan}
		}

		\begin{abstract}
			\NYY{
				We propose inverse problems of crack propagation using the phase-field
				models.
				First, we study the crack propagation in an inhomogeneous media
				in which fracture toughness varies in space.
				Using the two phase-field models based on different surface energy
				functionals, we perform simulations of
		 the crack propagation and show that the $J$-integral
		 reflects the effective inhomogeneous toughness.
				Then, we formulate regression problems to estimate
				space-dependent fracture toughness from the crack path.
				Our method successfully estimates the positions and magnitude of
				tougher regions.
				We also demonstrate that our method works for different geometry of inhomogeneity.
			}
		\end{abstract}
		
		
	\end{frontmatter}
	

	\section{Introduction}
	\NY{
		When materials are under strain, cracks propagate and the materials
		inevitably become fractured.
		Since the Griffith's theory \cite{Griffith_1920}, fracture mechanics
		has been developed.
		In this model, the critical toughness is interpreted as the balance
		between the elastic energy and surface energy of the material.
		Under the stain where the stored elastic energy is compatible with the
		surface energy, the material makes additional surface leading to the
		crack propagation.
		Although this picture needs to be extended to general crack phenomena,
		such as non-quasi-static crack propagation \cite{Karma:2004,RozenLevy_2020}, ductile fracture rather than
		brittle fracture \cite{Persson:2005,Ambati_2015b} and crack nucleation \cite{Tanne_2018}, the theory gives us a simple insight
		into the mechanism of crack propagation.
		
		One of the biggest difficulties of fracture mechanics is a singularity of
		stress at a crack tip.
		The solution of the linear elastic equation shows divergence $r^{-1/2}$ as the
		distance $r$ from the crack tip becomes smaller \cite{Anderson_2017}.
		This singularity suggests the presence of a microscopic length scale in
		the fracture phenomena.
		Several theoretical or numerical models have been proposed to handle the
		singularity, such as the extended finite element method (XFEM) \cite{Moes_1999}, the cohesive
		zone model \cite{Elices_2002}, and the phase-field model \cite{Wu_2019}.
		
		Among these methods, the phase-field method is particularly appealing
		because of its simple implementation of regularization parameter
		$\epsilon$ in the model.
		The phase-field model is based on the variational formula of the Griffith’s theoy \cite{Griffith_1920}.
		The theory states that a crack will propagate when the elastic energy released during the crack growth is greater than or
		equal to the surface energy which is proportional to the area of the new
		crack surfaces.
		Francfort and Marigo \cite{Francfort_Marigo_1998} proposed a variational
		method of crack
		propagation based on the minimization of the total energy, which is a sum of elastic
		energy and surface energy.
		In this approach, the bulk media and the crack surface are treated in a
		separate way; the former stores the elastic energy and the latter costs
		the surface energy.
		However, in numerical simulations, treatment of the discrete surface is
		difficult.
		In the phase-field model, the crack is expressed by the continuum
		variable $z({\bf x},t)$ in which $z \approx 1$ corresponds to crack
		regions whereas $z \approx 0$ corresponds to uncracked regions.
		Accordingly, both the elastic energy and surface energy are defined in
		the whole domain.
		The regularization parameter $\epsilon$ is introduced so that in the
		limit of $\epsilon \rightarrow 0$, the model converges to the
		Francfort-Marigo energy functional.
		
		Although they share the philosophy discussed above, there are several
		ways to formulate the phase-field models.
		In \cite{Bourdin_2000}, iterative minimization of the approximate Francfort-Marigo
		energy functional using the idea of Ambrosio-Tortorelli
		regularization was adapted.
		Crack propagation is realized by alternating minimization of the energy
		functional with respect to regularized crack field and strain field. 
		In \cite{Aranson:2000,Karma:2001,Hakim:2005}, crack propagation is
		described by a gradient flow of the Ginzuburg-Landau-type energy functional.
		Takaishi and Kimura \cite{Takaishi_Kimura_2009}, and Kuhn and
		M\"{u}ller \cite{Kuhn_2010} proposed a phase-field model which is
		remarkably simple for its treatment.
		It is described as a irreversible
		gradient flow of an approximate Francfort-Marigo energy with
		regularization.
		The irreversibility ensures that the crack does not heal.
		
		In this work, we study crack propagation in inhomogeneous materials.
		We consider the materials with tougher inclusions, and study how the
		crack propagation is modified by the inclusion, and how the effective
		toughness of the whole material is interpreted.
		To get deeper understanding of these issues, we propose the inverse
		problem of crack propagation in heterogeneous media.
		Crack propagation in inhomogenous materials have attracted a recent
		attention due to its relevance for understanding of fracture of
		composite materials.
		On the one hand, when the crack is propagating in a tougher region,
		the appearance of a new crack is suppressed and the material looks tougher.
		On the other hand, the crack may avoid the tougher region so that the
		crack does not feel that the material is tougher.
		Still, by the deformation of the crack path, the whole material may be
		considered tougher than the homogeneous one.
		In \cite{Bourdin_2014}, effective toughness was proposed for the
		materials with elastic inhomogeneity through the maximal $J$-integral
		during crack propagation.
		This idea was tested for the materials of epoxy resins \cite{Nishiura_2020}.
		The effective toughness of inhomogeneous toughness has been discussed
		theoretically and numerically in
		several studies \cite{Lebihain_Leblond_Ponson_2020}.
		
		Our aim in this study is to propose the inverse problem for crack
		propagation in inhomogeneous materials.
		We demonstrate the position and magnitude of tough regions can be
		estimated from the data of cracks.
		To formulate the inverse problem, we use the phase-field model based on
		the irreversible gradient flow.
		This method has an advantage that its regression formula can readily be
		obtained from the governing equations because the model is expressed by
		partial differential equations (PDEs).
		This is not the case for the method using iterative energy
		minimization.
		Within the approach, it is not clear how the inverse problem is formulated.

		Inverse problems of PDEs are inarguably
		important in a wide field of natural science.
		The goal is to estimate a model and parameters from data.
		When we have enough data of the left-hand side and right-hand side of
		our model, that is, variables and their time and spatial derivatives, we
		may use the regression to estimate parameters \cite{Baer:1999}.
		Using the sparse regularization, the model estimation, that is, the
		estimation of linear and nonlinear terms in the model to reproduce the
		data \cite{Rudy_Brunton_Proctor_Kutz_2017,Brunton_Kutz_2019}.
		Despite the success of the methods based on the regression for
		well-known nonlinear PDEs such as the Bugers equation and
		Kuramoto-Sivashinsky equation \cite{Rudy_Brunton_Proctor_Kutz_2017}, the
		application of the method to more complex PDEs remains challenging.
		In particular, it is known that the regression-based method is not
		robust against noisy data \cite{Schaeffer2017}.
		It also requires an accurate evaluation of time and spatial derivatives.
		In the problem of crack propagation, there is a sharp spatial change of
		the crack field $z({\bf x},t)$ to approximate the discrete crack, even
		if it is regularized.
		Another issue is irreversibility.
		Because of the plus operator in the model (see
		(\ref{Phase_field_model_TakaishiKimura}) and (\ref{Phase_field_model_Bourdin})), our model is
		not equality in the whole domain.
		Therefore, we need a pre-process to make a successful estimation of parameters.
		
		Crack propagation in heterogeneous media has less been studied using the
		irreversible gradient flow.
		Therefore, we also discuss forward problems.
		We consider several geometries of inhomogeneity, and how the crack path
		is dependent on them.
		We also study $J$-integral to clarify the effective toughness of each
		geometry.
		We consider two types of surface energy both for the
		forward and inverse problems.
		In \cite{Tanne_2018}, two types of surface energy both of which
		converge to Francfort-Marigo energy functional at $\epsilon \rightarrow
		0$ were discussed. 
		We will discuss the two surface energies may have advantages and
		disadvantages for the forward and inverse problems. 
	}

	The outline of this paper is as follows.
	\NYb{
		In Section
		\ref{section_problem_setting}, we present the two phase-field
		models, namely $\mathrm{AT}_1$ and $\mathrm{AT_2}$ model which are 
		based on the same elastic energy but two different types
		of surface energy.
		In Section \ref{section_crack_propagation}, we perform
		numerical simulations of the crack propagation in homogeneous
		and inhomogeneous media with different fracture toughness.
		We also compute the $J$-integral on the materials, and discuss
		how the $J$-integral depends on the crack path and reflects the
		effective fracture toughness. 
		In Section \ref{section_estimation}, we present the study of the
		estimation of the fracture toughness using the data of the crack
		path.
		We propose an algorithm based on data pre-processing, $k$-means
		classification, and linear regression.
		The method successfully estimates the position of the
		inhomogeneity of the fracture toughness and the position of
		tougher regions.
		Finally, in Section \ref{section_discussion}, we summarize the results and discuss the related topics.
	}	
	
	\section{Problem setting}\label{section_problem_setting}

	\NYb{
		We consider the crack growth phenomenon in space dimension two.
		We focus on the mode III fracture in which the
		anti-plane displacement is expressed by $u$.
		Following the spirit of the phase-field models, the crack is
		described by the order parameter $z$ such that
		$z \in[0,1]$,
		\begin{equation*}
			\begin{cases}
				z \approx 1 \quad \quad  \mbox{cracked}\\
				z \approx 0 \quad \quad  \mbox{not cracked.} 
			\end{cases}
		\end{equation*}
		We consider a bounded domain $\Omega=(0,L)\times(0,H)$ in space dimension two, whose
		boundary $\Gamma$ is piecewise smooth. Let $\Gamma_{\mathrm{D}}$
		be an open subset of $\Gamma$ with a finite number of connected
		componenents, and we define $\Gamma_{\mathrm{N}}=\Gamma
		\backslash \Gamma_{\mathrm{D}}$.
		The phase-field models are based on the variational formula in
		which we minimize the total energy $\mathcal{E} = \mathcal{E}_1
		+ \mathcal{E}_2$ consisting of elastic and
		surface energy.
		The regularized elastic energy by $\displaystyle{\mathcal{E}_1}$, and it is defined as 
		\begin{equation}\label{def_elastic_energy}
			\begin{aligned}
				&\mathcal{E}_1(u,z) \\
				:= &	\displaystyle{\frac{\mu}{2} \int_{\Omega} (1-z)^2 \lvert \nabla u \rvert^2 \mathrm{d}x - \int_{\Omega}fu\mathrm{d}x -\int_{\Gamma_{\mathrm{N}}} (1-z)^2hu\mathrm{d}s             } 
			\end{aligned}
		\end{equation}
		where  $\mu$ is the elastic constant and is one of the Lam\'e
		constant.
		Here, $f$ and $h$ are external forces.
		In this work, we do not consider these external forces, and therefore, $f=h=0$.
		
		The two types of surface energy with different power of $z$
		 are discussed in \cite{Tanne_2018}.
		The two models are referred to as $\mathrm{AT}_1$ and
		$\mathrm{AT_2}$, and expressed as
		\begin{equation}\label{surface_energy_AT1}
			\displaystyle{\mathcal{E}_2(z)  =\frac{3}{8} \int_{\Omega} \gamma(\mathbf{x}) \left( \epsilon \lvert \nabla z\rvert^2   +\frac{1}{\epsilon}z   \right)  \mathrm{d}\mathbf{x} } \quad (\mathrm{AT_1})
		\end{equation}
		and
		\begin{equation}\label{surface_energy_AT2}
			\displaystyle{\mathcal{E}_2(z) = \frac{1}{2} \int_{\Omega} \gamma(\mathbf{x}) \left( \epsilon \lvert \nabla z\rvert^2   +\frac{1}{\epsilon}z^2   \right)  \mathrm{d}\mathbf{x} } \quad (\mathrm{AT_2}).
		\end{equation}
		In Takaishi-Kimura \cite{Takaishi_Kimura_2009} and Kuhn and
		M\"uller \cite{Kuhn_2010}, the authors consider the
		$\mathrm{AT_2}$ model and derive the gradient flow of
		$\mathcal{E}_1+\mathcal{E}_2$ and then apply the plus operator
		to ensure the propagation of the crack.
		Then, the model is given by 
		\begin{equation}\label{Phase_field_model_TakaishiKimura}
			\begin{cases}
				\displaystyle{{\alpha_1 \frac{\partial u}{\partial t} }
					= \mu \operatorname{div}\left( (1-z)^2 \nabla u \right)}\\
				\displaystyle{\alpha_2 \frac{\partial z}{\partial t}} =
				{\Bigg(}\epsilon \operatorname{div}( \gamma(\mathbf{x})\nabla z) 
				- \frac{\gamma(\mathbf{x})}{\epsilon}z+ \mu \lvert \nabla u \rvert^2 (1-z)  {\Bigg)_{+}}
			\end{cases}
		\end{equation}
		with boundary conditions $u(\mathbf{x},t) = g(\mathbf{x},t)$ on
		$ \Gamma_{\mathrm{D}} \times (0,T)$,
		$\displaystyle{\frac{\partial u}{\partial n} =0}$ on $
		\Gamma_{\mathrm{N}} \times (0,T)$  for $u$.
		Here, $g$ expresses the boundary displacement.
		We also apply homogeneous Neumann boundary conditions
		$\displaystyle{\frac{\partial z}{\partial n} =0}$ on $\Gamma
		\times (0,T)$ for $z$.
		Here, $ \Gamma_{\mathrm{D}} =(0,L)\times \{x_2=0,x_2=H\}$ and $ \Gamma_{\mathrm{N}} = \{x_1=0, x_1= L\}\times (0,H)$.
		We denote spatial position as $\mathbf{x} = (x_1,x_2) \in
		\mathbb{R}^2$.
		The plus operator $(a)_+=\max(a,0)$ acts on the gradient of the
		total energy with respect to $z$.
		This operator is introduced to express non-repairing crack, and thus, $\displaystyle{\frac{\partial z}{\partial t}=0}$ when the total energy is increasing.
			%
		
		By applying the same idea of deriving the gradient flow of the
		total energy and then applying the plus operator $()_+$, we may propose the following system which corresponds to the $\mathrm{AT_1}$ surface energy:
		\begin{equation}\label{Phase_field_model_Bourdin}
			\begin{cases}
				\displaystyle{{\alpha_1 \frac{\partial u}{\partial t} }
					= \mu \operatorname{div}\left( (1-z)^2 \nabla u \right)} \quad & \\
				\displaystyle{\alpha_2 \frac{\partial z}{\partial t} =
					{\Bigg(} \frac{3}{4}\epsilon \operatorname{div}( \gamma(\mathbf{x})\nabla z) 
					- \frac{3}{8}\frac{\gamma(\mathbf{x})}{\epsilon} 
					+ \mu \lvert \nabla u \rvert^2 (1-z)  {\Bigg)_{+}}},
			\end{cases}
		\end{equation}
		with the same boundary settings as mentioned in \eqref{Phase_field_model_TakaishiKimura}.
		In the following, we call the two models as the $\mathrm{AT_1}$
		model \eqref{Phase_field_model_Bourdin} and as the $\mathrm{AT_2}$
		model \eqref{Phase_field_model_TakaishiKimura}.
		We will first present results based of the $\mathrm{AT_2}$ model
		and then the results of the $\mathrm{AT_1}$ model. 
		
		Both in the two phase-field models, the parameters $\alpha_1$ and
		$\alpha_2$ describe dissipative relaxation time scales.
		Because we focus on quasi-static crack propagation, we consider
		the stationary elastic equation, namely $\alpha_1=0$.
		We also set the time scale for the crack field as $\alpha_1=10^{-3}$.
		In the spirit of the phase-field models, $\epsilon$ is a
		regularization parameter such that the minimum length scale of
		$z$ is given as $O(\epsilon)$.
		We apply the surfing boundary conditions \cite{Bourdin_2014, Nishiura_2020} for the two models, that is 
		\NYYY{
			\begin{equation}\label{BC_strain}
				g(\mathbf{x},t)=\frac{{A}}{2}\left(1-\tanh \left(\frac{x_1-{v}t}{{d}}\right)\right)\operatorname{sign}(x_2-0.5H).  
			\end{equation}}
		In the following, we consider the space domain $(0,5) \times (0,1)$, namely, $L=5$ and $H=1$. The parameter values in the model as well as in the boundary
		condition will be specified in the next section. 
	}

	
	%
	\section{Crack Propagation in heterogeneous
		media}\label{section_crack_propagation}
	\NY{
		In this section, we consider crack propagation of a material with
		tougher inclusions.
		The tougher regions have a larger value of $\gamma$ so that more
		elastic energy should be stored to break this region.}	In the following, we denote by $\gamma_0$ the fracture toughness of the homogeneous media, whereas $\gamma_1$ for the fracture toughness of the inclusions.
	\NY{
		We will discuss how the crack is propagating in the heterogeneous media, and
		how to characterize the effective toughness of the heterogeneous
		material using $J$-integral.
		Before showing the numerical results, we the discuss qualitative
		classification of crack paths.
		
		In \cite{Lebihain_Leblond_Ponson_2020}, Lebihain et al. classified every point along the crack front in four different states.
		Inspired by this idea, we classify crack propagation into
		four categories: (i) propagation in homogeneous material, (ii)
		penetration, (iii) being stuck, and (iv) bypassing.
		In (i), the crack simply does not hit an inclusion, and move in
		a homogeneous region. 
		In (ii)-(iv),  the crack interacts with the inhomogeneity.
		When a crack reaches an inclusion, the crack moves inside the
		inclusion in (ii).
		In (iii), the crack stops either inside the inclusion or at the interface
		between homogeneous media and the inclusion.
		In (iv), the crack avoids the inclusion.
		In the following numerical simulations, we will see by applying different geometry of inclusion and different toughness of the inclusion, we observe the crack path in the four categories.
		We observe (i) straight propagation in homogeneous material;
		(ii) penetration in inhomogeneous material with periodic stripes
		or one-disk inclusion with $\gamma_1$ not large;
		(iii) being stuck in periodic stripes or one-disk inclusion with
		$\gamma_1$ large;
		and (iv) bypassing in one-disk inclusion case with $\gamma_1$ large.

	}

	\subsection{Simulations of crack propagation}\label{subsection_simulation_path}

	
		\begin{figure}[htbp] 
		\begin{center}
			\includegraphics[scale=0.3]{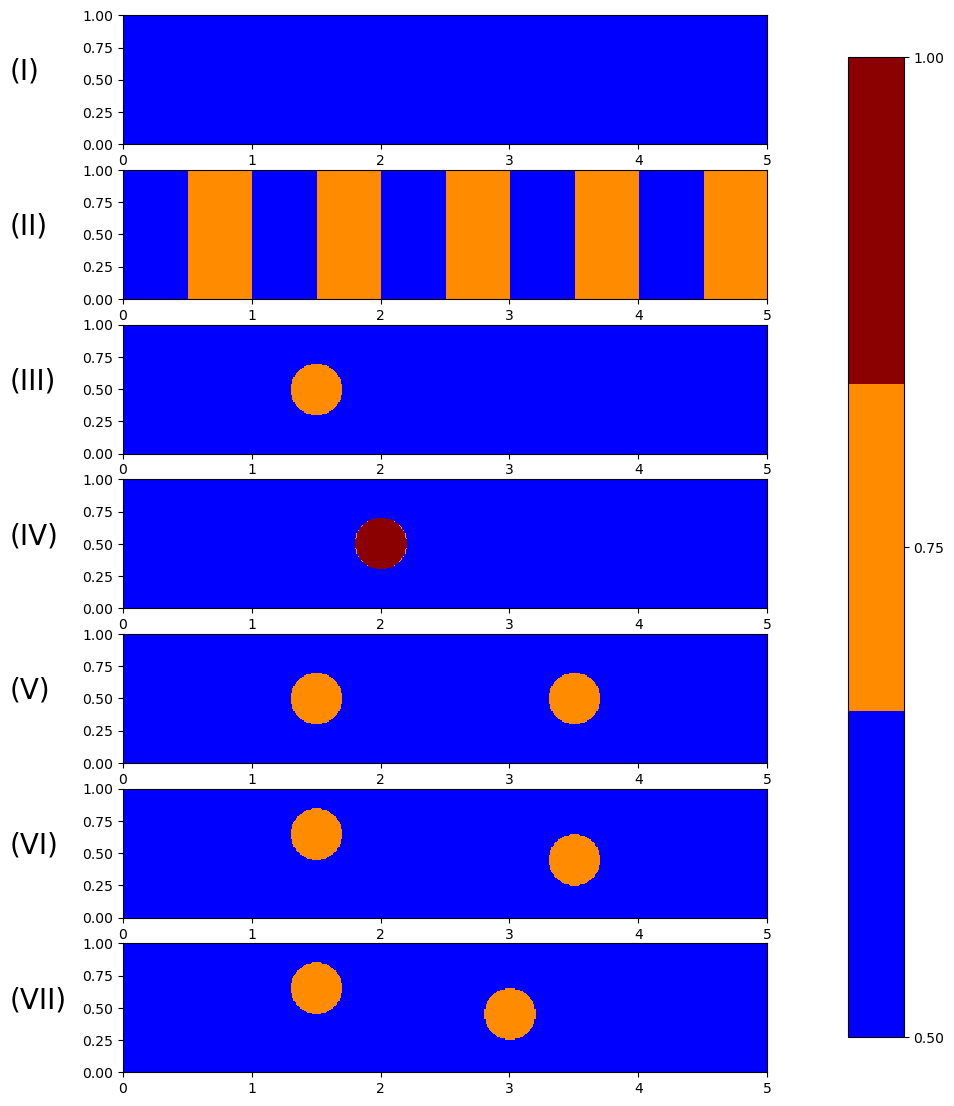}
		\end{center}
		\caption{
			\NYY{
				The different distribution of inhomogeneous fracture toughness
				that we consider. (I) homogeneous media, (II) stripe inclusion,
				(III-VII) disk inclusions.
				The color indicates fracture toughness; homogeneous media is
				shown in blue, whereas tougher inclusions are shown in red.
				One disk inclusions with (III) $\gamma_1=0.75$ and (IV)
				$\gamma_1=1.0$ are considered.
				Two disk inclusions with different positions in $x_1$ and $x_2$ are shown
				in (V-VII). 
			}
		}
		\label{Fracture_Toughness_6cases}
	\end{figure}
	
	\NY{
		In order to study the difference among different geometry of
		inclusion and different toughness of the inclusion, we first
		consider the following four cases, namely (I) homogeneous toughness with
		\NYY{$\gamma_0=0.5$}; (II) periodic stripes whose toughness is
		\NYY{$\gamma_1=0.75$} and the width of each region is 0.5;
		(III) one-disk inclusion with a disk whose barycenter is at $(1.5,0.5)$
		with radius 0.2 and toughness $\gamma_1=0.75$; and (IV) one harder disk inclusion
		with disk's barycenter at $(2.0,0.5)$ with radius 0.2 and
		\NYY{$\gamma_1=1.0$}.
		For the inhomogeneous materials, we set \NYY{$\gamma_0=0.5$} outside the
		inclusion, which we call homogeneous media.
		We show the inhomogeneous distribution of $\gamma (\mathbf{x})$ in Fig.~\ref{Fracture_Toughness_6cases}.
		
			We should note that stripe inclusion and disk inclusion are
		inherently different because bypassing occurs only for the disk inclusion.		
		We also remark that both in the simulations and in the estimation of fracture toughness, we neglect the interfacial effect of the inclusion, which was included in Lebihain et al. \cite{Lebihain_Leblond_Ponson_2020}.
		
		In addition to the four cases,  we also consider three cases of
		two-disks inclusions with different positions.
		The inclusion disks' coordinates are (V) $(1.5,0.5),(3.5,0.5)$;
		(VI) $(1.5,0.65),(3.5,0.45)$ and (VII) $(1.5,0.65),(3,0.45)$ with
		$\gamma_1=0.75$.
		\NYY{Figure~\ref{Fracture_Toughness_6cases}}(V)-(VII) shows $\gamma ({\bf x})$ for
		the three cases.
		These examples indicate that in (V) the two inclusions are
		placed at the center in the $x_2$ direction, and therefore the
		system has symmetry around the $x_1$ axis.
		In (VI) and (VII), the $x_2$ positions of the two inclusions are
		shifted from the center.
		Because the initial crack is chosen at the center of the $x_2$
		axis and the crack goes straight in the homogeneous media, the
		crack reaches different positions on the surface of the
		inclusions.
		Therefore, the angle between the incoming crack to the
		inclusion, and line connecting the center of the inclusion and
		the its surface which the crack reaches, is different
		for each inclusion in (VI) and (VII).
		When the crack bypasses the inclusion, the
		crack path is expected to be deformed.
		The incoming angle of a crack reaching the second (right) inclusion may be
		changed by the distance between the two inclusions.
		We demonstrate it by using (VI) and (VII).
		
	}
		\NYb{
		We perform numerical simulations of the phase-field models
		\eqref{Phase_field_model_TakaishiKimura} and \eqref{Phase_field_model_Bourdin} by a finite volume
		method.
		We refer to \cite{Eymard_Herbin_Gallouet_2010} for the method as well as for the notation of the mesh.	
		The discrete solution of $z$ and $u$ are denoted by $\{Z_K^n\}$ and $\{U_K^n\}$ over control volume $K$ in time interval $[t_n, t_{n+1})$ respectively. We apply the uniform space mesh, which is the square size of $0.01 \times 0.01$ and 
		we apply uniform time discretization, that is we fix
		the time step to be $\Delta t = 10^{-2}$ with $t_n = n \Delta t$ so that $[0,T)= \cup_{n=0}^{N_T-1}[
		t_n,t_{n+1})$. 
		We set the parameters as $\epsilon = 0.02$ and $\mu=1$. 
		For the surfing boundary condition, we first consider $v= 1$ and
		$d=0.5$.
		We will specify the applied strain $A$ in each case.
	}

	\subsubsection*{Crack path and tip position in the
		$\mathrm{AT_2}$ model}
	
\begin{figure}[htbp]
		\includegraphics[scale=0.3]{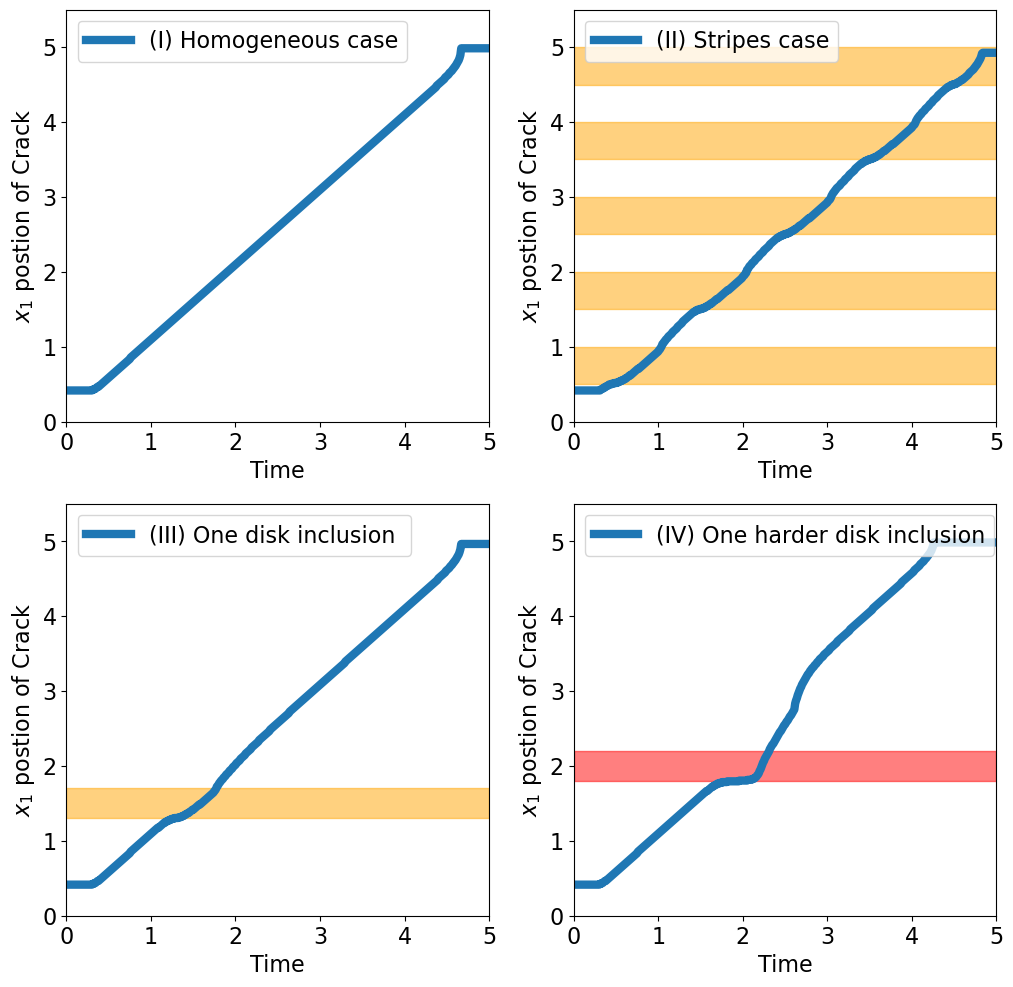}
		\includegraphics[scale=0.25]{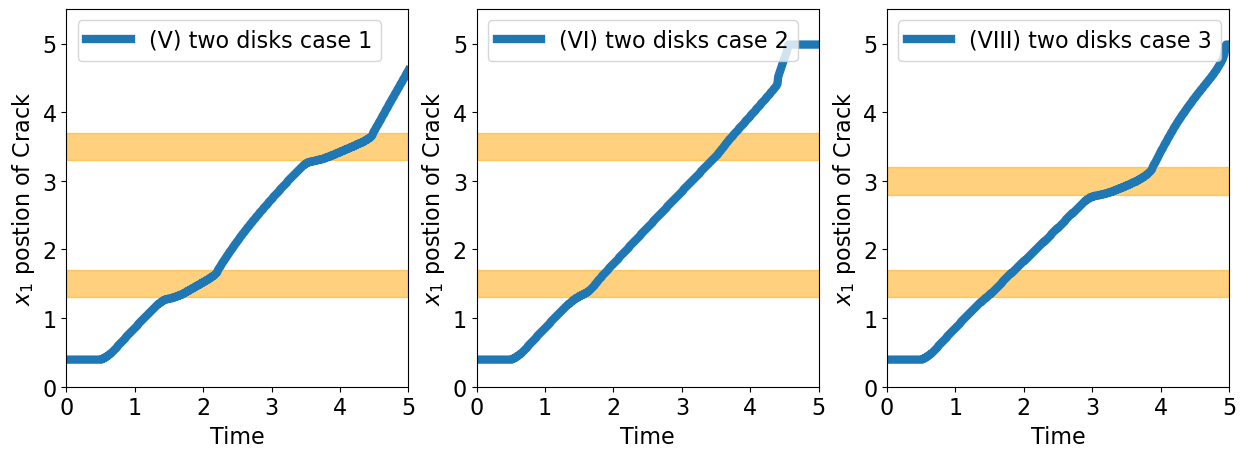}
		\caption{
			\NYY{
				Positions of crack tips in $x_1$ as a function of time using
				the $\mathrm{AT_2}$ model.
				The shaded area corresponds to the inhomogeneity in
				Fig.~\ref{Fracture_Toughness_6cases}.
				\label{fig_crack_tip_at2}
			}
		}
	\end{figure}

	\NY{
		We present the simulation results of the crack paths in
		the $\mathrm{AT_2}$ model, described in
		\eqref{Phase_field_model_TakaishiKimura} with the initial condition as mentioned in \cite{Takaishi_Kimura_2009}. 
		The crack paths in the cases (I)-(IV) are shown in
		Fig.\ref{Fig_Allpaths_AT2}.
		The crack paths are described by $z({\bf x})$.
		In the cases (I)-(III), the paths are straight, and therefore,
		penetration occurs for (II) and (III).
		On the other hand, the crack avoids in (IV), and it corresponds
		to bypassing.			
		\NYY{Figure~\ref{Fig_Allpaths_AT2}} shows the crack paths in the
		cases (V)-(VII).
		In (V), the crack penetrates the two inclusions, whereas the
		crack avoids two inclusions in (VI).
		In (VII), the crack penetrates the right inclusion, but bypasses
		the left inclusion.
		It is clear that the shapes of crack paths are dependent both on
		the distribution of tougher regions, magnitude of the toughness,
		and their positions.
		
		\begin{figure}[htbp] 
			\begin{center}
				\includegraphics[scale=0.3]{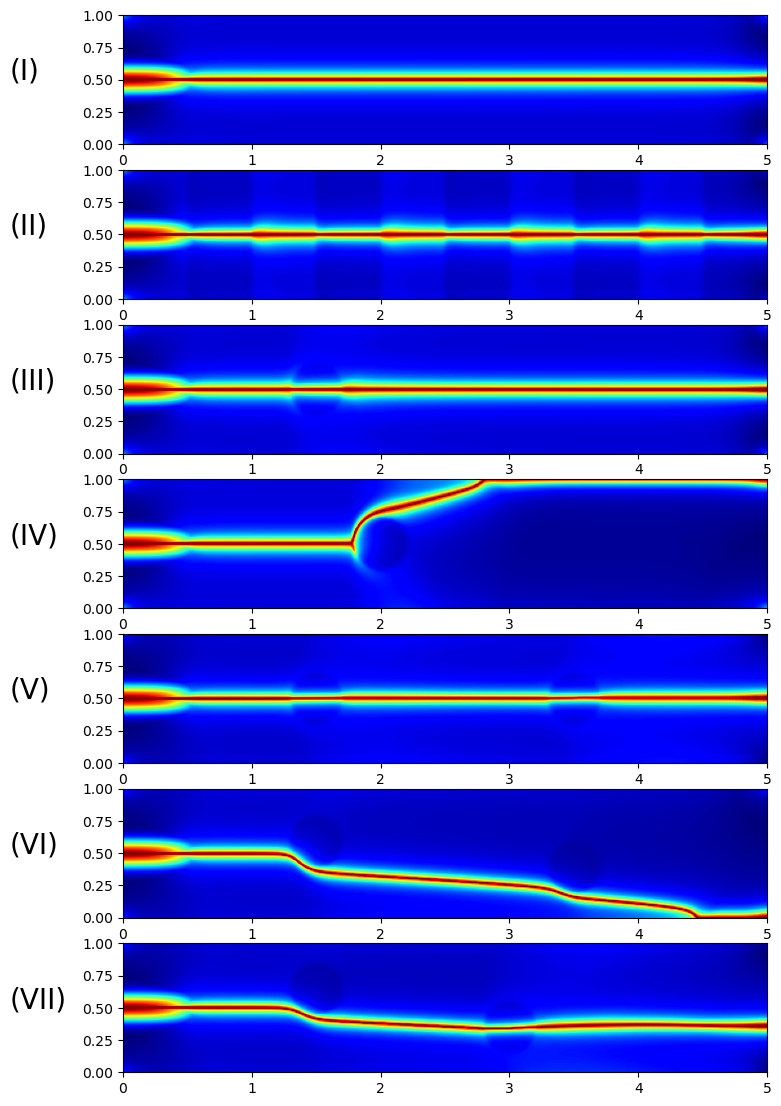}
			\end{center}
			\caption{
				\NYY{
					Crack paths in inhomogeneous materials with different spatial
					dependence of the fracture toughness $\gamma$ using the
					$\mathrm{AT_2}$ model.
					These paths correspond to the inhomogeneity in Fig.~\ref{Fracture_Toughness_6cases}.
				}
				\NYYY{
					The applied strain $A$ in \eqref{BC_strain} is chosen as $A=1.25$ for
					(I)-(IV), and $A=1$ for (V)-(VII).
				}
			} \label{Fig_Allpaths_AT2}
		\end{figure}

		As we will see later, the position of a crack tip, $x_1$, is an important observable for the $J$-integral computation and for the inverse problem, we therefore present the crack tip position in different cases. 
		\NYY{Figure~\ref{fig_crack_tip_at2}} shows the position of a crack tip $x_1$ as a function of time.
		When the crack propagates in the homogeneous media, the speed of the
		propagation is approximately constant, and therefore, $x_1$ increases
		linearly in time.
		When the crack reaches an inclusion, the speed of the propagation
		decreases.
		After penetrating the inclusion, the speed of crack increases, and it
		becomes even higher than that of the crack in the homogeneous media.
		When the crack reaches a harder inclusion in \NYY{Fig.~\ref{Fig_Allpaths_AT2}(IV)}, the crack
		propagation stops, and when bypassing, the speed of the crack
		increases \NYY{(see Fig.~\ref{fig_crack_tip_at2}(IV))}.
		
		For the systems with two inclusions, the position of a crack tip
		$x_1(t)$ shows a similar behavior to the case of one inclusion.
		The speed of crack propagation slows down when it moves inside
		the inclusion, and the speed increases after penetration.
		In contrast with \NYY{Fig.~\ref{Fig_Allpaths_AT2}(IV)}, when the
		crack bypasses the inclusion as in
		\NYY{Fig.~\ref{Fig_Allpaths_AT2}(VI) and (VII)}, the speed does
		not change when the crack reaches the inclusion
		\NYY{(Fig.~\ref{fig_crack_tip_at2}(VI) and (VII))}.
		This is because the incoming angle of the crack deviates from
		$\theta=0$.
		In this case, bypassing is more likely to occur, and
		accordingly, the crack does not get stuck when it reaches inclusion.
	}
	
	\subsubsection*{Crack path and tip position in the
		$\mathrm{AT_1}$ model}
		
\begin{figure}[htbp] 
			\begin{center}
				\includegraphics[scale=0.3]{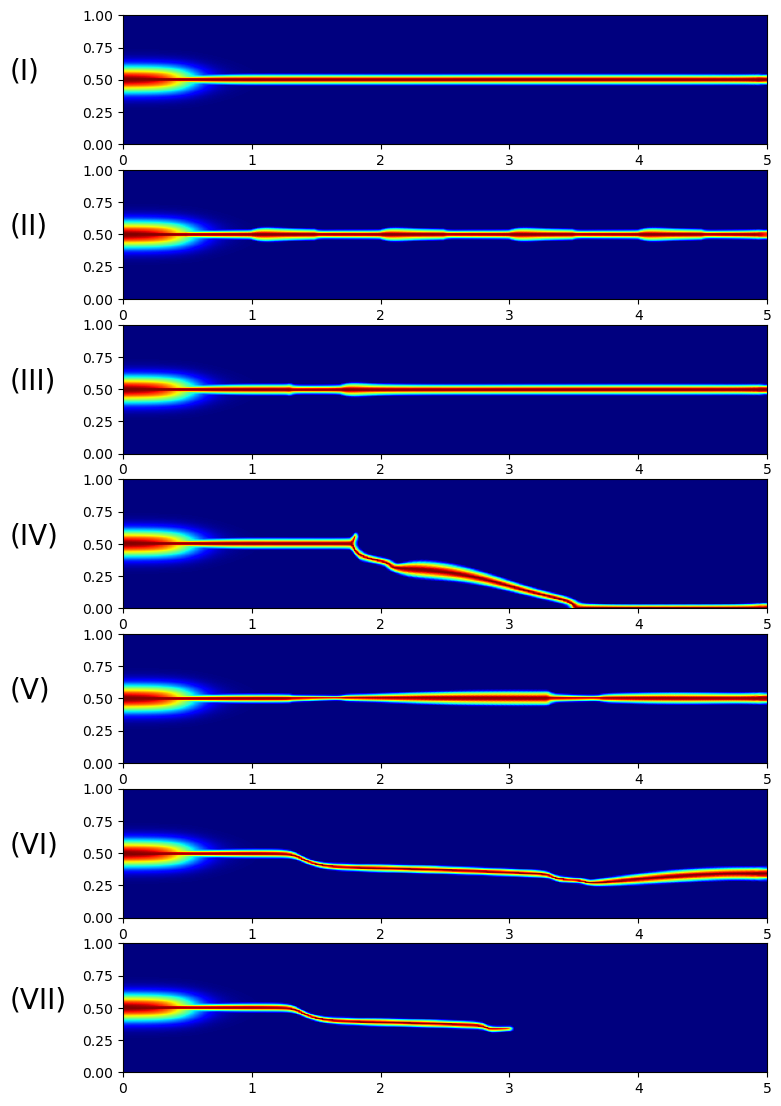}
			\end{center}
			\caption{
				\NYY{
					Crack paths in inhomogeneous materials with different spatial
					dependence of the fracture toughness $\gamma$ using the
					$\mathrm{AT_1}$ model.
					These paths correspond to the inhomogeneity in Fig.~\ref{Fracture_Toughness_6cases}.
				}
				\NYYY{
					The applied strain $A$ in \eqref{BC_strain} is chosen as $A=1.25$ for
					(I)-(IV), and $A=1$ for (V)-(VII).
				}		 
			} \label{Fig_crackpath_AT1_First4cases}
		\end{figure}

\begin{figure}[htbp]
			\includegraphics[scale=0.3]{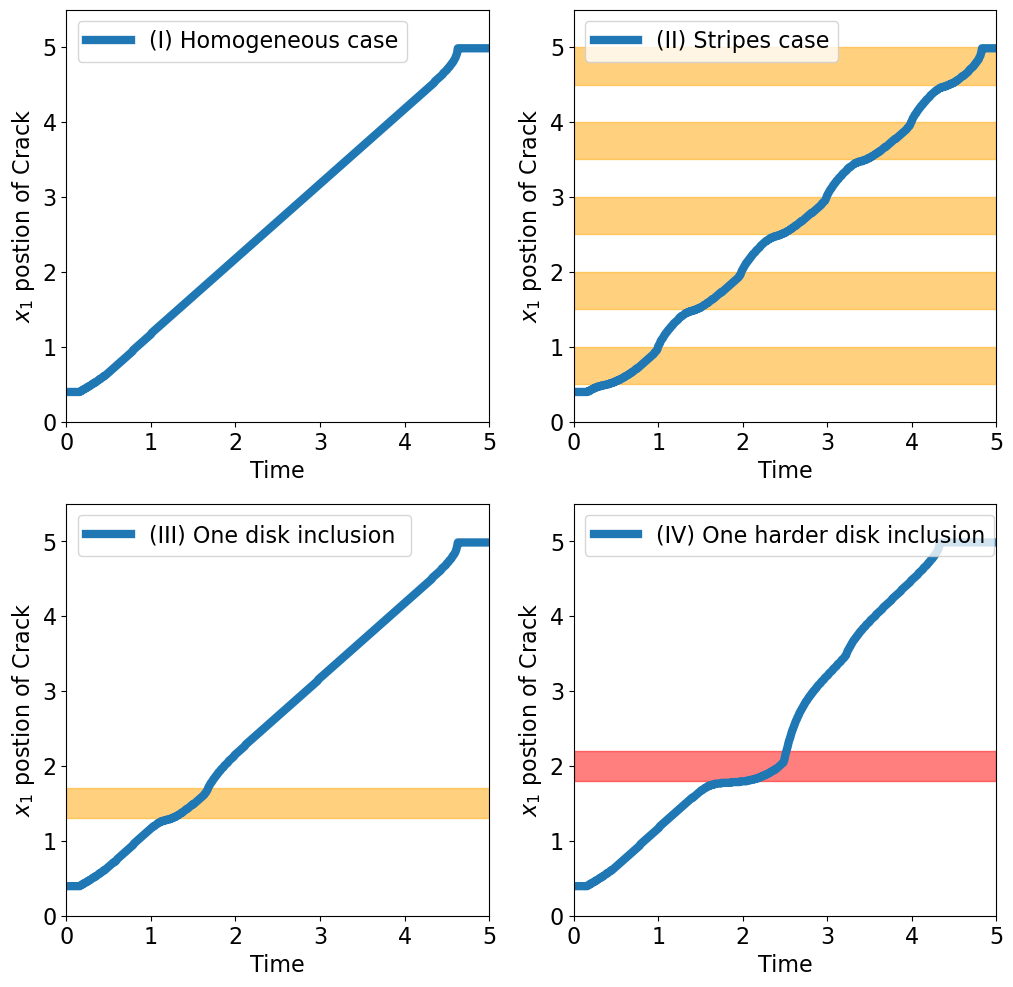}
			\includegraphics[scale=0.25]{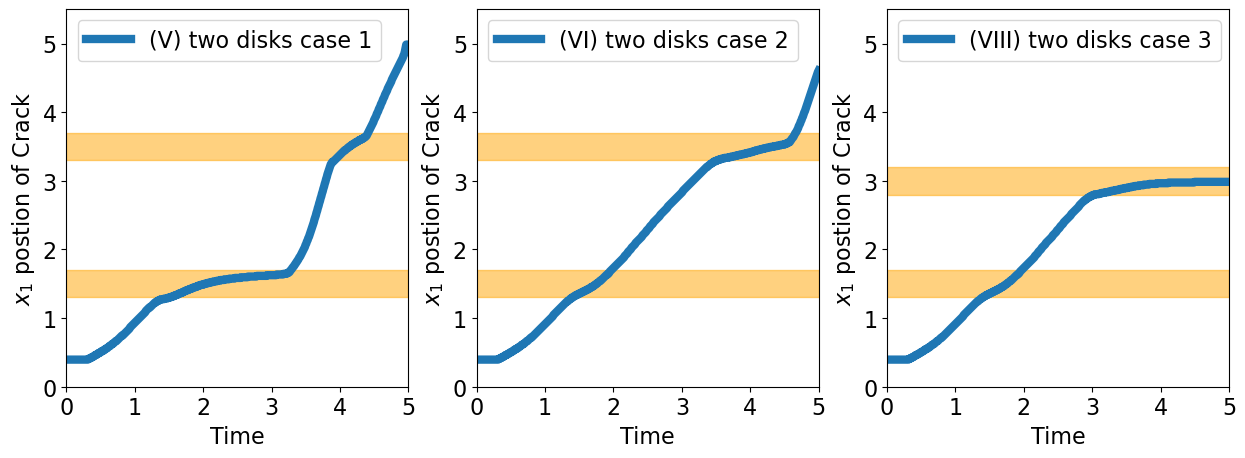}
			\caption{
				\NYY{
					Positions of crack tips in $x_1$ as a function of time using
					the $\mathrm{AT_2}$ model.
					The shaded area corresponds to the inhomogeneity in
					Fig.~\ref{Fracture_Toughness_6cases}.
					\label{fig_crack_tip_at1}
				}
			}
		\end{figure}

		
		\NY{
			Next, we consider the crack propagation in a heterogeneous media in the $\mathrm{AT_1}$ model.
			We use the same inhomogeneity as shown in Fig.~\ref{Fracture_Toughness_6cases} 
			with the same surfing boundary conditions.
			The crack paths obtained by this model are qualitatively
			the same as those obtained by the $\mathrm{AT_1}$ model.
			The results are shown in \NYY{Fig.\ref{Fig_crackpath_AT1_First4cases}}.
			In the case of (VII), the crack is eventually stuck at
			the right inclusion, and correspondingly $x_1$ saturates
			\NYY{(see Fig.~\ref{fig_crack_tip_at1})}.
		}
		
				%
		
\begin{figure}[htbp]
	\includegraphics[scale=0.25]{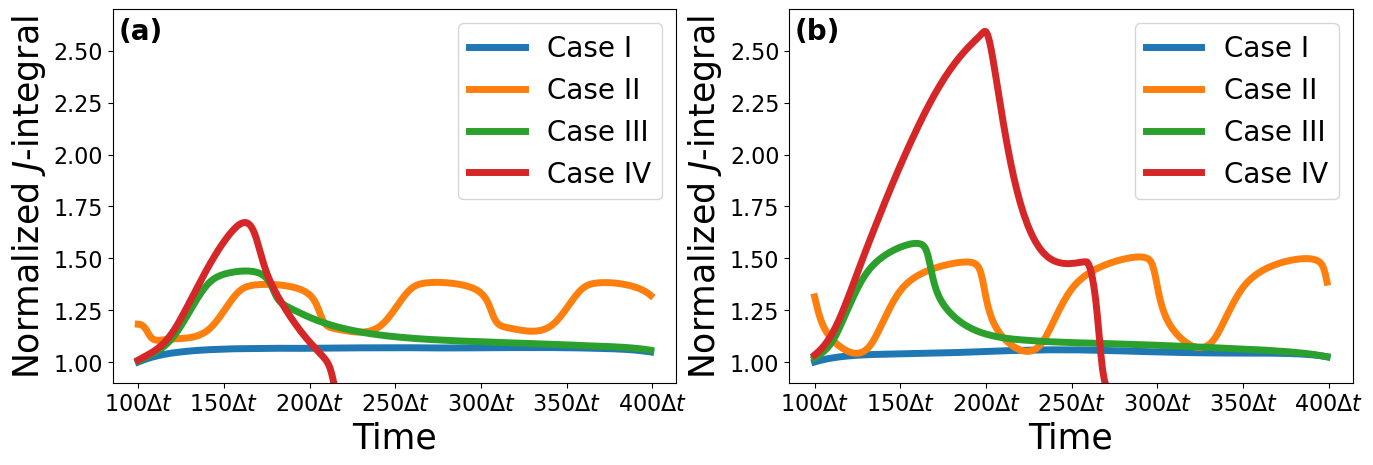}
			\caption{
				\NYY{The $J$-integral during crack propagation in inhomogeneous media
					using the (A) $\mathrm{AT_2}$ and
					(B) $\mathrm{AT_1}$ models
		 corresponding to Fig.~\ref{Fracture_Toughness_6cases}(I)-(IV).
				}
				\NYYY{
					The applied strain $A$ in \eqref{BC_strain} is chosen as $A=1.25$.
				}		 
			} \label{J_integral_AT2}
		\end{figure}

		\NY{
			\subsection{characterization of effective toughness by
				$J$-integral during crack propagation}\label{subsection_J_integral}
			
			In the previous section, we observe different crack
			paths depending on the toughness of the inclusion, and
			its shape and geometry.
			In this section, we discuss
			effective toughness as a global property of the toughness of the
			heterogeneous material.
			For this purpose, we consider the $J$-integral.
			The $J$-integral is the strain energy release rate defined as
			\begin{equation*}\label{def_J_integral}
				\begin{aligned}
					J :&= - \frac{d (\mbox{Potential Energy)}}{d (\mbox{Crack Area})}\\ &=\int_{{\Gamma_J}}\left({W}dy-{\mathbf{T}}\cdot \frac{\partial \mathbf{u}}{\partial x_1}ds \right),
				\end{aligned}
			\end{equation*}
			where $\Gamma_J$ is an arbitrary {counter}-clockwise curve
			around the crack tip as shown in Fig.~\ref{J_integral_boundary}, $W$ is the strain energy density,
			and $\mathbf{T} = (T_1,T_2,T_3)$ is the surface traction vector.
			For a homogeneous media, the $J$-integral corresponds to surface energy
			density $\gamma$ in
			\eqref{surface_energy_AT1}
			and \eqref{surface_energy_AT2}.
			In this case, the $J$-integral is
			independent of the choice of the
			integration path $\Gamma_J$ \cite{Anderson_2017}.
			On the other hand, in the inhomogeneous material, the $J$-integral is path-dependent.
			
			\NYb{
				We show the computation of the $J$-integral in the
				four cases (I)-(IV) that we consider in
				the Fig.~\ref{Fracture_Toughness_6cases}.
				We choose the integration path $\Gamma_J$ in a
				counter-clockwise manner as in
				Fig.~\ref{J_integral_boundary}.
				Then, the $J$-integral is decomposed by the integral values on each boundary  such that $J=J_{\Gamma_1}+J_{\Gamma_2}+J_{\Gamma_3}+J_{\Gamma_4}$. 			
			}
			We refer to \cite{Nishiura_2020, Bourdin_2014} for the
			numerical computation of the $J$-integral both in
			homogeneous and inhomogenous systems using the
			phase-field model.
			
			A direct computation yields 
			\begin{equation}\label{J_Gamma1_Formula}
				J_{\Gamma_1}= - \int_{\Gamma_1}(1-z)^2\left( \frac{\partial u}{\partial x_2}\right)^2dx_2,
			\end{equation}
			and
			\begin{equation}\label{J_Gamma3_Formula}
				J_{\Gamma_3} = \int_{\Gamma_3}(1-z)^2 \frac{\partial u}{\partial x_2}\frac{\partial u}{\partial x_1}ds,
			\end{equation}
			and because of the symmetry,
			$J_{\Gamma_4}=J_{\Gamma_3}$.
			Because of the surfing boundary condition on 
			$\Gamma_3$ and on $\Gamma_4$, when a crack propagates towards $\Gamma_2$, no strain will be applied
			on $\Gamma_2$, so that $J_{\Gamma_2}$ can be
			approximately ignored.
			When the crack approaches $\Gamma_2$, $J_{\Gamma_2}$ becomes nonzero, 
			indicating that the finite system size affects the
			computation of the $J$-integral.
			Therefore, in the following, we present the computation of the
			$J$-integral before the moment $J_{\Gamma_2}$ starts to increase.
			
			\begin{figure}[htbp] 
				\begin{center}
					\includegraphics[scale=0.3]{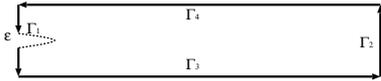}
				\end{center}
				\caption{
					\NYY{
						The contours $\Gamma_J$ for the computation of
						the $J$-integral.
						The cracked region is shown by the curve in the dotted line.
						The width of the region is $\sim \epsilon$ in
						\eqref{Phase_field_model_TakaishiKimura} and \eqref{Phase_field_model_Bourdin}.
					}
				}
				\label{J_integral_boundary}
			\end{figure}

			The value of $J$-integral is sensitive to the choice of
			the position of $\Gamma_1$ because of
			the large value of ${\partial
				u}/{\partial y}$.
			We do not set $\Gamma_1$ to be the boundary $\{x_1=0\} \times (0,H)$ of the domain because the initial condition do not satisfy the equation.
			In the following, both in $\mathrm{AT_2}$ and  $\mathrm{AT_1}$ models, we set the position of $\Gamma_1$ to be $49 \Delta x$.

			\NYb{
				Figure~\ref{J_integral_AT2} shows the
				$J$-integral during propagation.
				We normalize the $J$-integral by the value for
				the homogeneous media at $t=100 \Delta t$ to see
				the dependence on the toughness $\gamma$ in the
				inclusion.
			}
			For the homogeneous media ((I) in Fig.\ref{Fracture_Toughness_6cases}), the $J$-integral is constant in time.
			This result is consistent with the constant speed of
			crack propagation.
			When the crack propagates in tough inclusions ((II)-(IV) in
			Fig.\ref{Fracture_Toughness_6cases})), the $J$-integral
			increases.
			After going out from the inclusions, the $J$-integral
			decreases to the value of the homogeneous media.
			In the case of (II) and (III), the maximum value is
			approximately $\gamma$ of the tougher regions,
			{$\gamma_1/\gamma_0 \approx 1.5$}.
			When the toughness of the inclusion is larger, the
			maximum value of the $J$-integral is larger.
		}
		
		
		%
		

		%
		
\begin{figure}[htbp]
			\centering
			\includegraphics[scale=0.25]{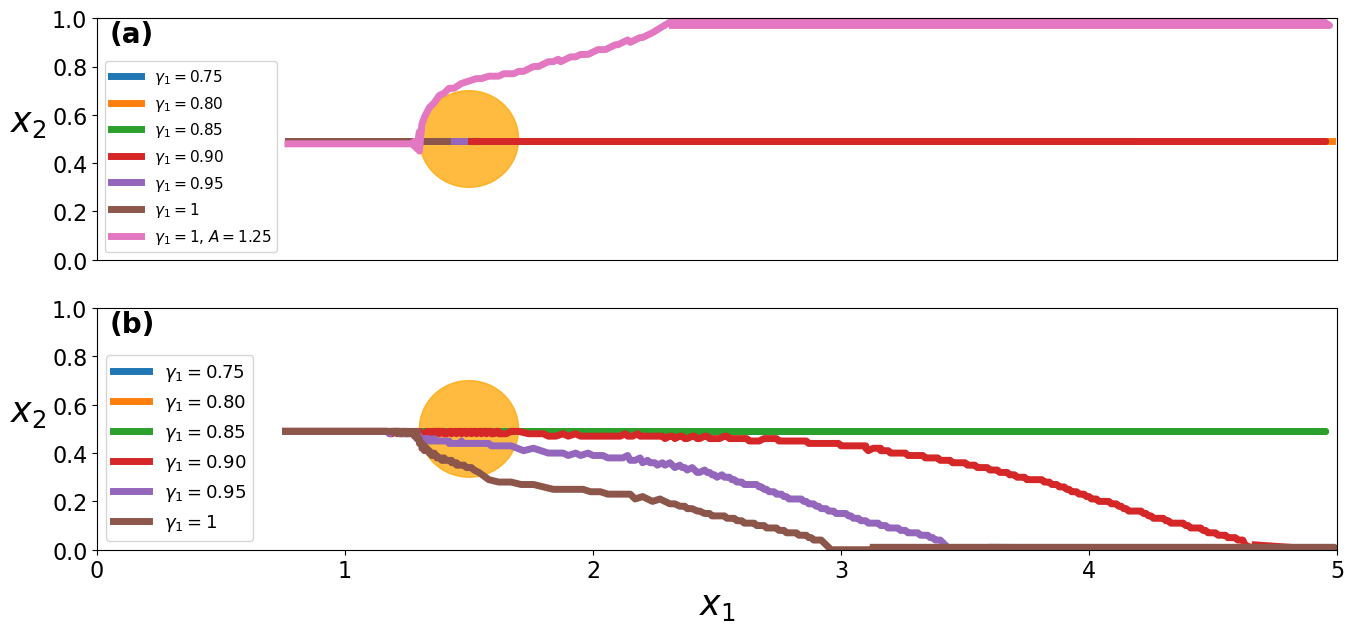}
			\caption{
				\NYY{
					The crack paths obtained by the tip positions
					$(x_1,x_2)$ in the inhomogeneous media under different
					toughness $\gamma_1$.
					The results of the (a) $\mathrm{AT_2}$
					and (b) $\mathrm{AT_1}$ models are shown. {The applied strain $A$ of \eqref{BC_strain} is chosen as $A=1$ in (a) and $A=1.25$ in (b).}
					\label{Fig_Jintegral_gamma1}
				}
			}
		\end{figure}

		\NYYY{
				To study how inhomogeneous toughness affects crack propagation, we
			focus on the cases of the one-disk inclusion whose toughness is gradually
			increased from $0.75$ to $1.0$.
			We start from $\mathrm{AT_2}$ model with a disk inclusion at $(1.5,0.5)$ with radius
			0.2.
			The toughness of the inclusion is $0.75, 0.8, 0.85, 0.9,
			0.95$ to $1$ with the $A$ in the surfing boundary
			condition to be $1$.
			We find that the crack path changes from penetration to
			pinning as the toughness of the inclusion increases, as
			shown in Fig.~\ref{Fig_Jintegral_gamma1}.
			When $\gamma = 1$ and $A=1.25$, the crack path shows
			bypassing after pinning at the surface of the inclusion. 
			\GY{We then consider $\mathrm{AT_1}$ model with $A=1.25$. 
			We find that the crack paths change from straight penetration to deviated penetration as the toughness of the inclusion increases.}
			
			The results of corresponding $J$-integrals to each path
			in Fig.~\ref{J_integral_AT2} are shown in Fig.\ref{Fig_Jintegral_gamma1}.
			In both models, the peak value of $J$-integral increases
			as $\gamma_1$ in the disk inclusion increases.
			Before the $J$-integral reaches the value comparable with
			$\gamma_1/\gamma_0$, the crack is pinned at the surface
			of the inclusion.
			When $\gamma_1/\gamma_0$ is small, the crack starts to move after
			reaching the peak, and results in penetration.
			On the other hand, when $\gamma_1/\gamma_0$ is large,
			the $J$-integral saturates and the crack remains at the
			pinned state.
			When the applied strain is larger $A=1.25$ in
			Fig.\ref{Fig_Jintegral_gamma1}(a), the $J$-integral
			increases as the crack tip hits the disk and is pinned before bypassing occurs.
			After the crack starts to move.
			The $J$-integral reaches its peak and takes a larger
			value than $\gamma_1/\gamma_0$.
			After the crack starts to move, the $J$-integral drops.
			The crack path is no longer along the $x_1$ direction,
			and therefore, the $J$-integral does not come back to
			the value of the inhomogeneous media.
			These behaviors are shared both by the $\mathrm{AT_1}$ and $\mathrm{AT_2}$ models.
		}

		\begin{figure}[htbp]
			\includegraphics[scale=0.25]{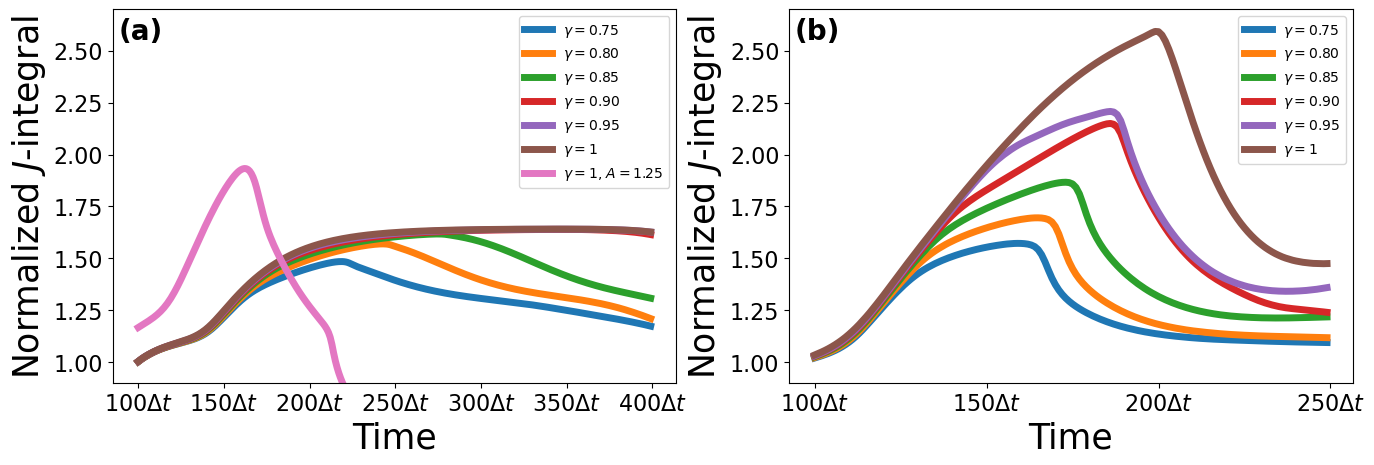}
			\caption{
				\NYY{
					The  $J$-integral during crack propagation in the inhomogeneous media under different
					toughness $\gamma_1$ corresponding to Fig.~\ref{Fig_Jintegral_gamma1}.
					The results of the (a) $\mathrm{AT_2}$
										and (b) $\mathrm{AT_1}$ models are shown. {The applied strain $A$ of \eqref{BC_strain} is chosen as $A=1$ in (a) and $A=1.25$ in (b).}
				}
			}
		\end{figure}

\begin{figure}[htbp]
			\begin{center}
				\includegraphics[scale=0.25]{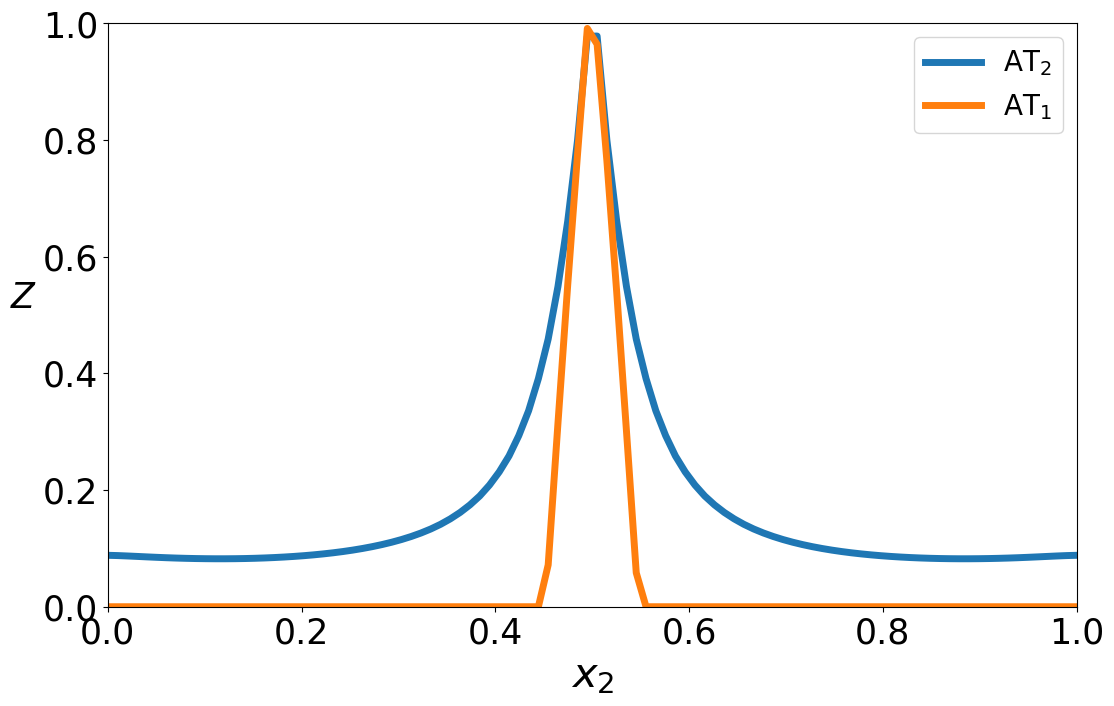}
			\end{center}
			\caption{
				\NYY{
					The cross sections along $x_2$ of the crack field $z(x_2)$ in
					the homogeneous media using the $\mathrm{AT_2}$ and
					$\mathrm{AT_1}$ models at time $t=174\Delta t$.
				}
			} \label{Solution_Cut_AT2}
		\end{figure}

		\NYb{
			Despite the similar results of the $\mathrm{AT_1}$ and $\mathrm{AT_2}$
			models, there are differences between these models.		
			In Fig.~\ref{Solution_Cut_AT2}, we check the shape of
			the solutions. 
			Given a time $t^n$, we obtain the simulation result of
			$z(x_1,x_2)$, then fix the $x_1$ coordinates and plot
			the solution of $z$ as a function of $x_2$ coordinates.
			We find that for the solution in the $\mathrm{AT_2}$
			model, the values on $\{x_2 = 0, x_2 = 1 \}$, which
			correspond to the boundary $\Gamma_3$ and $\Gamma_4$,
			are not $0$.
			The difference was also discussed in Tanne et al. \cite{Tanne_2018} and Wu et
			al. \cite{Wu_2019}.
			In Wu et al. \cite{Wu_2019}, the authors show the distribution
			of the crack phase-field in different models.
			In our notation, the solution is approximated, for fixed $x_1$, as
			$z=(1-\frac{|x_2|}{2l_0})^2$ in  $\mathrm{AT_1}$, whereas
			$z=\exp(-\frac{|x_2|}{l_0})$ in $\mathrm{AT_2}$, where $l_0$ is
			the length scale. This is also shown in our case in Figure
			\ref{Solution_Cut_AT2}.
			Accordingly, the effect of the plus operator is very different
			in the two models, that is in $\mathrm{AT_1}$, the effect is
			over the whole domain, whereas in $\mathrm{AT_2}$, the effect is
			mainly on the region where the value of $z$ is large.		
			This difference may affect the inverse estimation of the fracture toughness that we present in Sec.~\ref{section_estimation}.				
		}
		
		\section{Estimation of fracture toughness}\label{section_estimation}
		
		\NY{
			The aim of this section is to propose and demonstrate
			the method to estimate the fracture toughness $\gamma$
			from the crack path obtained from the phase-field model. 
			More precisely, we perform simulation for the
			phase-field models, and apply a regression technique on
			numerical results $\{Z_K^n\}$ and $\{U_K^n\}$ to
			estimate $\gamma$.
			The basic idea is to optimise $\gamma$ so that the
			left-hand and right-hand sides of
			\eqref{Phase_field_model_TakaishiKimura} or
			\eqref{Phase_field_model_Bourdin} become closer in a
			certain norm.
			This is not straightforward as it looks because
			of the nonlinearlity caused by the $(\cdot)_+$ operator.
			In fact, when the $(\cdot)_+$ operator works,
			\eqref{Phase_field_model_TakaishiKimura} and
			\eqref{Phase_field_model_Bourdin} without $(\cdot)_+$
			operator are not equality but inequality.
			Therefore, we cannot use all the sample points in the
			system, but need pre-process to choose appropriate
			sample points.
			
			\subsection{inverse problems of the $\mathrm{AT_2}$ model}\label{section_Inverse_problem_Takaishi_Kimura}
			We start from the model of $\mathrm{AT_2}$,  and we propose the following algorithm.
			\begin{enumerate}
				\item \textbf{Simulation}:
				We perform the simulation of the phase-field model of $\mathrm{AT_2}$.
				At the same time, we track the position of the crack tip
				at each time step.
				We denote the position of the crack tip as
				$(x_{1,n},x_{2,n} )$ at time step $n\Delta t$ for
				all $1 \leq n \leq N_T$.		       
				\item \textbf{Sampling}:
				At each time step, we
				sample the grid points in space {in the region
					where the crack has not arrived}.
				We consider the time interval $[n_0\Delta t, n_1
				\Delta t)$.
				We define $\Omega_n = [x_{1,n},L]\times [0,H]$,
				that is the space region where the crack has not
				arrived yet.
				And then we sample $N_{x_1}$ grid points over
				$[x_{1,n},L]$ and $N_{x_2}$ grid points on
				$[0,H]$.
				The $j^{th}$ sample point will be denoted by $p_j:=((x_{1}^{(j)},x_{2}^{(j)}),t^{(j)})$ and the total sample set is denoted by $P =\{p_j\}_{j=1,2,...,N_{data}}$ with $N_{data} = (n_1-n_0)\cdot N_{x_1} \cdot N_{x_2}$. 
				
				\item \textbf{Interpolation}:
				We apply Chebyshev polynomial of degree $7$ for
				the interpolation of $\{Z(P)\}$ and $\{ U(P)\}$.
				The results of the interpolation is denoted by
				$\widetilde{Z}(x_1,x_2,t)$ and
				$\widetilde{U}(x_1,x_2,t)$.
				Then, by computing {the derivative of the polynomial}, we obtain the values of $\displaystyle{{\partial \widetilde{Z}}/{\partial t}}$, $\displaystyle{{\partial \widetilde{Z}}/{\partial {x_1}}}$, $\displaystyle{{\partial \widetilde{Z}}/{\partial {x_2}}}$, $\Delta \widetilde{Z}$ and $ \lvert \nabla \widetilde{U}\rvert^2 (1-\widetilde{Z})$.
				
				\item \textbf{Linear Regression}:
				We assume that the values of $\alpha_2$, $\epsilon$, and $\mu$ are known and that $\gamma$ is constant, which corresponds to the homogeneous material. 
				Based on the second equation of system
				\eqref{Phase_field_model_TakaishiKimura},  if we
				omit the plus operator $(\cdot)_+$, formally we
				obtain
				\begin{equation}\label{def_estimate_TakaishiKimura}
					\alpha_2  \frac{\partial
						\widetilde{Z}}{\partial t}  =\gamma
					\left(\epsilon \Delta \widetilde{Z} -
					\frac{\widetilde{Z}}{\epsilon}\right) +\mu
					\lvert \nabla
					\widetilde{U}\rvert^2(1-\widetilde{Z}),
				\end{equation}
				from which we define
				\begin{equation}\label{def_Y}
					Y:= \alpha_2  \frac{\partial \widetilde{Z}}{\partial t} - \mu \lvert \nabla \widetilde{U}\rvert^2 (1-\widetilde{Z}) 
				\end{equation}
				and 
				\begin{equation}\label{def_X}
					X :=  \epsilon \Delta \widetilde{Z}-\frac{\widetilde{Z}}{\epsilon}.
				\end{equation}
				
				\textcolor{black}{If we find the linear relation such that $Y=
					\hat{\gamma} X $, then
					$\hat{\gamma}$ is the estimator of $\gamma$, that
					is the fracture toughness of the phase-field
					model.
					We will generalize this approach to the
					inhomogeneous $\gamma$.}
		\end{enumerate}}
		
\begin{figure}[htbp] 
				\begin{center}
					\includegraphics[scale=0.3]{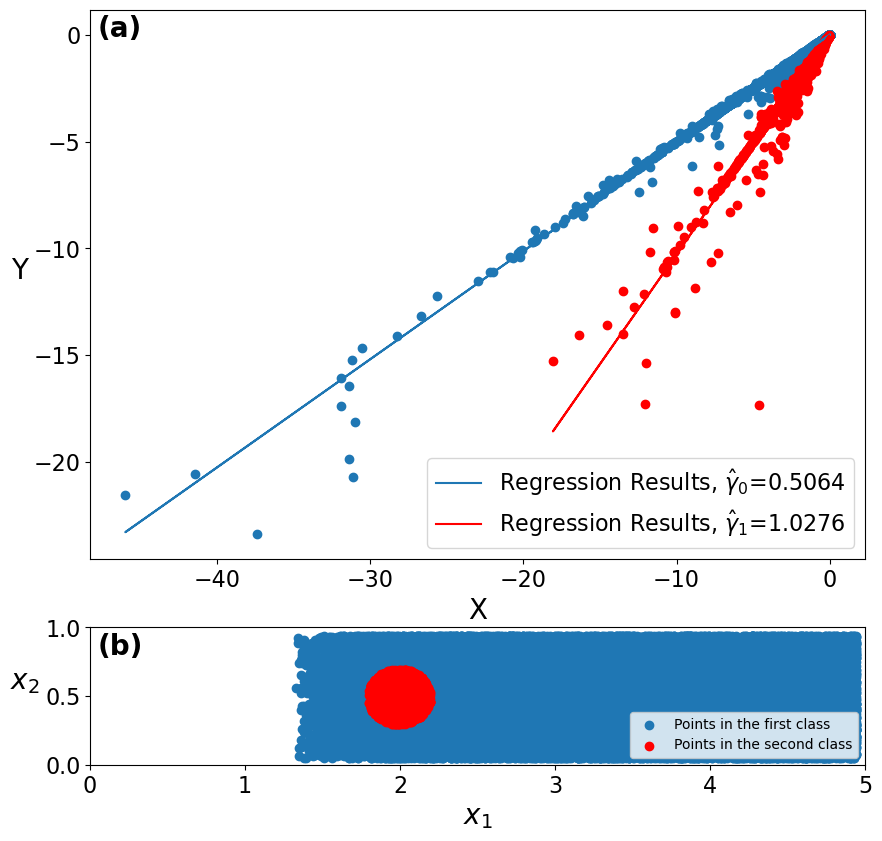}
				\end{center}
				\caption{
					\NYY{
						(a) The sampled data in $(X,Y)$ obtained by the
						$\mathrm{AT_2}$ model for the crack
						propagation in inhomogeneous media with a
						harder disk
						inclusion.
						The data is classified by using the $k$-means algorithm.
						The solid line is the result of
						linear regression for each class.
						(b) The positions of the data point
						corresponding to each class. 
					}
					\NYYY{
						The applied strain $A$ in \eqref{BC_strain} is chosen as $A=1.25$.
					}		 		 
				}  \label{1harderdisk_gamma_difclasses_AT2}
			\end{figure}

		
		\NY{	\subsubsection{homogeneous case}
			
			We start from the numerical results of a homogeneous material with $\gamma_0 =0.5$.
			We consider the time interval $[100\Delta t, 200\Delta t)$ and at each time step we sample $20\times 20$ the
			grid points.
			The positions of sample points are shown in
			Fig.~\ref{Total_sample_homo_treated}.
			Based on the sample points, we perform interpolation and
			compute time and spatial derivatives of $Z$ to evaluate
			$X$ and $Y$.
			The results are shown in \NYY{Fig.~\ref{Total_sample_homo_treated}}.
			Applying the linear regression, we obtain $Y=0.5002
			\cdot X $.
			The ground truth of $\gamma_0$ is $\gamma_0=0.5$, and thus we obtain an accurate estimation.		
			\begin{figure}[htbp] 
				\begin{center}
					\includegraphics[scale=0.3]{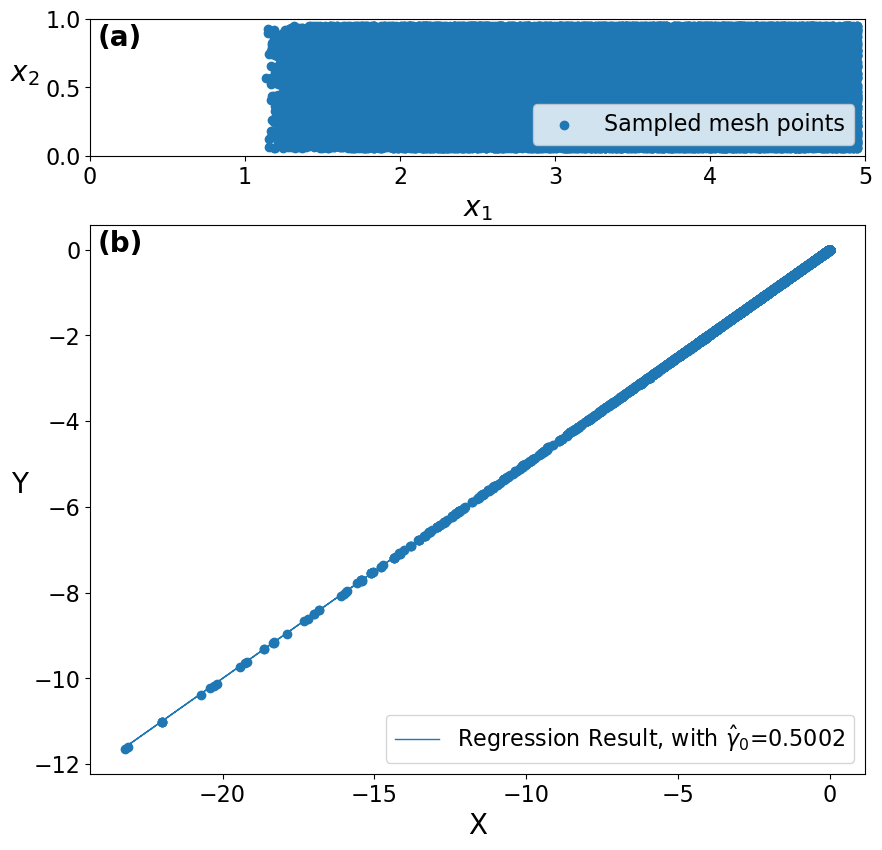}
				\end{center}
				\caption{
					\NYY{
						(a) Sampled grid points at which crack
						has not yet arrived.
						(b) The sampled data points in $(X,Y)$ of
						\eqref{def_X} and \eqref{def_Y}.
						The solid line is the result of linear regression.
					}
					\NYYY{
						The applied strain $A$ in \eqref{BC_strain} is chosen as $A=1.25$.
					}		 			 
				} \label{Total_sample_homo_treated}
			\end{figure}

			\subsubsection{Inhomogeneous case}\label{section_Inverse_problem_Bourdin}
			
			In the inhomogeneous case, the fracture toughness
			$\gamma (x)$ is non-uniform in space.
			Therefore, we expand $\operatorname{div}(\gamma(x)\nabla
			z)$ and omit the plus operator $(\cdot)_+$ in
			\eqref{Phase_field_model_TakaishiKimura} to obtain
			\begin{equation}
				\begin{aligned}
					\alpha_2  \frac{\partial \widetilde{Z}}{\partial t}  =&\gamma (x)\left(\epsilon \Delta \widetilde{Z} - \frac{\widetilde{Z}}{\epsilon}\right) 
					+\mu \lvert \nabla \widetilde{U}\rvert^2(1-\widetilde{Z})\\
					&+ \epsilon \frac{\partial \gamma(x)}{\partial {x_1} } \cdot \frac{\partial \widetilde{Z}}{\partial {x_1} } +   \epsilon\frac{\partial \gamma(x)}{\partial {x_2} }\cdot \frac{\partial \widetilde{Z}}{\partial {x_2} }.
				\end{aligned}
			\end{equation}		
			Similar to step 4 of the homogeneous case, we define 
			\begin{equation*}
				Y:= \alpha_2  \frac{\partial \widetilde{Z}}{\partial t} - \mu \lvert \nabla \widetilde{U}\rvert^2 (1-\widetilde{Z})
			\end{equation*}
			and
			\begin{equation*}
				X:=  \epsilon \Delta \widetilde{Z}-\frac{\widetilde{Z}}{\epsilon},\quad X_1:= \epsilon \frac{\partial \widetilde{Z}}{\partial {x_1} },\quad X_2:=  \epsilon\frac{\partial \widetilde{Z}}{\partial {x_2} }.
			\end{equation*}
			Because we consider the inclusion that has larger value of $\gamma$,
			which is sharply varied from the value in the media, the
			terms ${\partial \gamma(x)}/{\partial {x_1} }$ and
			${\partial \gamma(x)}/{\partial {x_2} }$ 
			approximately vanish at the many
			sample points.
			Therefore, we may apply the regression using
			$Y$ and $X$ to analyse the inhomogeneity in the same way as we did in the homogeneous case.
			
			We assume that the fracture toughness $\gamma$ can take
			two values; \NYY{$\gamma_0$} in the media and \NYY{$\gamma_1$} in
			an inclusion.
			Although each sample point of $(X,Y)$ does not know
			whether this point is \NYY{$\gamma_0$ or $\gamma_1$}, we may
			classify each point using the method of modified
			$k$-means (see \ref{section_kmeans} for
			the detailed method).
			As in the homogeneous case, the grid points are sampled in the region where the crack has not arrived in each case.
			
			\NYb{
				Similar to Sec.\ref{section_crack_propagation},
				we consider the inhomogeneity of case (II)-(VII)
				in Fig.~\ref{Fracture_Toughness_6cases}.
				We choose the time interval to be $[100\Delta t,
				200\Delta t)$, except the case (IV) in which we consider the time interval $[120\Delta t, 200\Delta
				t)$ because the position of the inclusion is $x_1=2.0$
				(\NYYY{see
					Fig.~\ref{1harderdisk_gamma_difclasses_AT2}}).
			}
			
			\begin{figure}[H] 
				\begin{center}
					\includegraphics[scale=0.2]{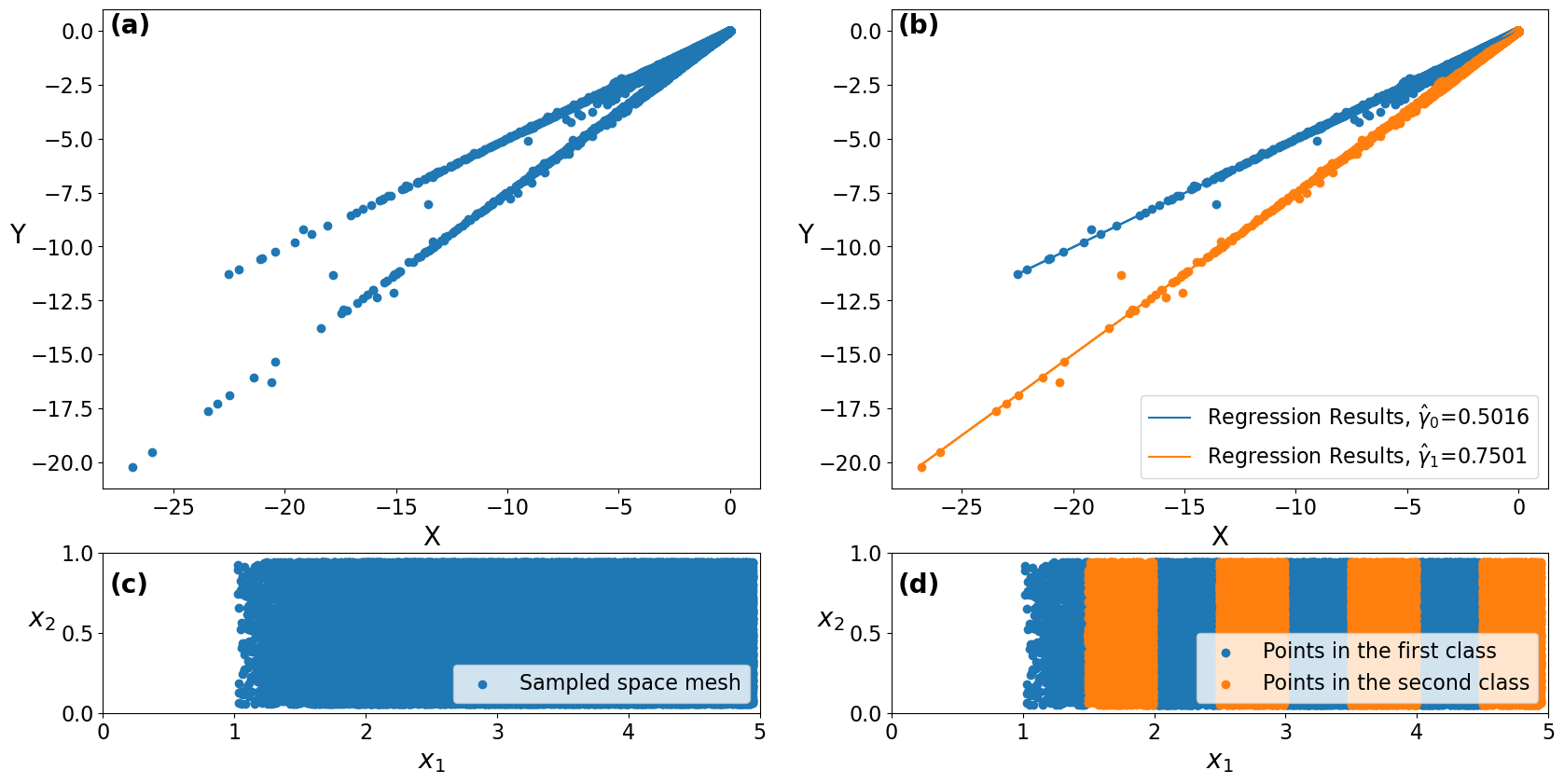}
				\end{center}
				\caption{
					\NYY{
						(a,b) The sampled data in $(X,Y)$ for the crack
						propagation in inhomogeneous media with stripe
						inclusions before (a) and after (b) the
						classification using the $k$-means algorithm.
						The solid line in (b) is the result of
						linear regression for each class.
						(c,d) The positions of the data point in real
						space before (a) and after (b) the classification. 
					}
					\NYYY{
						The applied strain $A$ in \eqref{BC_strain} is chosen as $A=1.25$.
					}		 			 
				}  \label{Stripes_points}
			\end{figure}
			
			\NYb{			
				We first consider the stripe inclusions.
				We show in Fig.~\ref{Stripes_points}(a) the results of $X$ and $Y$ after sampling and interpolation.
				The result indicates that there are two slopes.
				In fact, by applying $k$-means, we may classify the points and estimate the $\gamma$ in each class.
				The estimated values of $\gamma$ are \NYY{$\hat{\gamma}_0=0.5016$
					and $\hat{\gamma}_1=0.7501$} which are close to the
				ground truth, $0.5$ and $0.75$.
				In Fig.\ref{Stripes_points}, we show the spatial
				position of the points in the two classes.
				The result shows a good agreement with the spatial
				distribution of $\gamma(x)$ in
				Fig.\ref{Fracture_Toughness_6cases}(II).
				This result suggests that the method can discover not only the different values of the inhomogeneity, but also the position of the inhomogeneity.

				Next, we consider the one-disk inclusion.
				In Fig.\ref{1disk_2diskscase1_difclasses_AT2}(a,b), we show the
				results of the sampling after interpolation and
				$k$-means method are applied.
				The estimated values of $\gamma$ are
				\NYY{$\hat{\gamma}_0=0.5008$ and
				$\hat{\gamma}_1=0.7534$} compared to
				$0.5$ and $0.75$ of the ground truth.
				Our method also works for the disk inclusion with harder
				toughness
				(Fig.\ref{Fracture_Toughness_6cases}(IV)).			
				The estimated values of $\gamma$ are
				\NYY{
					$\hat{\gamma}_0=0.5064$ and $\hat{\gamma}_1=1.0276$ compared to their
					ground truth values $0.5$ and
					$1.0$, respectively.}
				The result is shown in
				Fig.~\ref{1harderdisk_gamma_difclasses_AT2} in Supplementary Materials.
			}

			Our method to estimate the toughness $\gamma$ of
			inclusions and their positions also work for multiple inclusions.
			We consider the three cases (V-VII) in \NYY{Fig.\ref{Fracture_Toughness_6cases}}.
			\NYb{Figure~\ref{1disk_2diskscase1_difclasses_AT2}(c,d)} show the results of
			the case (V)
			after the pre-processing of data, interpolation and $k$-means method as mentioned in the
			algorithm.
			The results show that we could accurately estimate two
			$\gamma$ corresponding to the toughness of media
			\NYY{, $\gamma_0$, and inclusions, $\gamma_1$}.
			The estimated toughness is \NYY{$\hat{\gamma_0}=0.4963$ and
				$\hat{\gamma_1}=0.7562$.}
			Note that in this case, the crack penetrates both
			inclusions (see \NYY{Fig.~\ref{Fig_Allpaths_AT2})}.
			
			\begin{figure}[htbp] 
				\begin{centering}
					
					\includegraphics[scale=0.2]{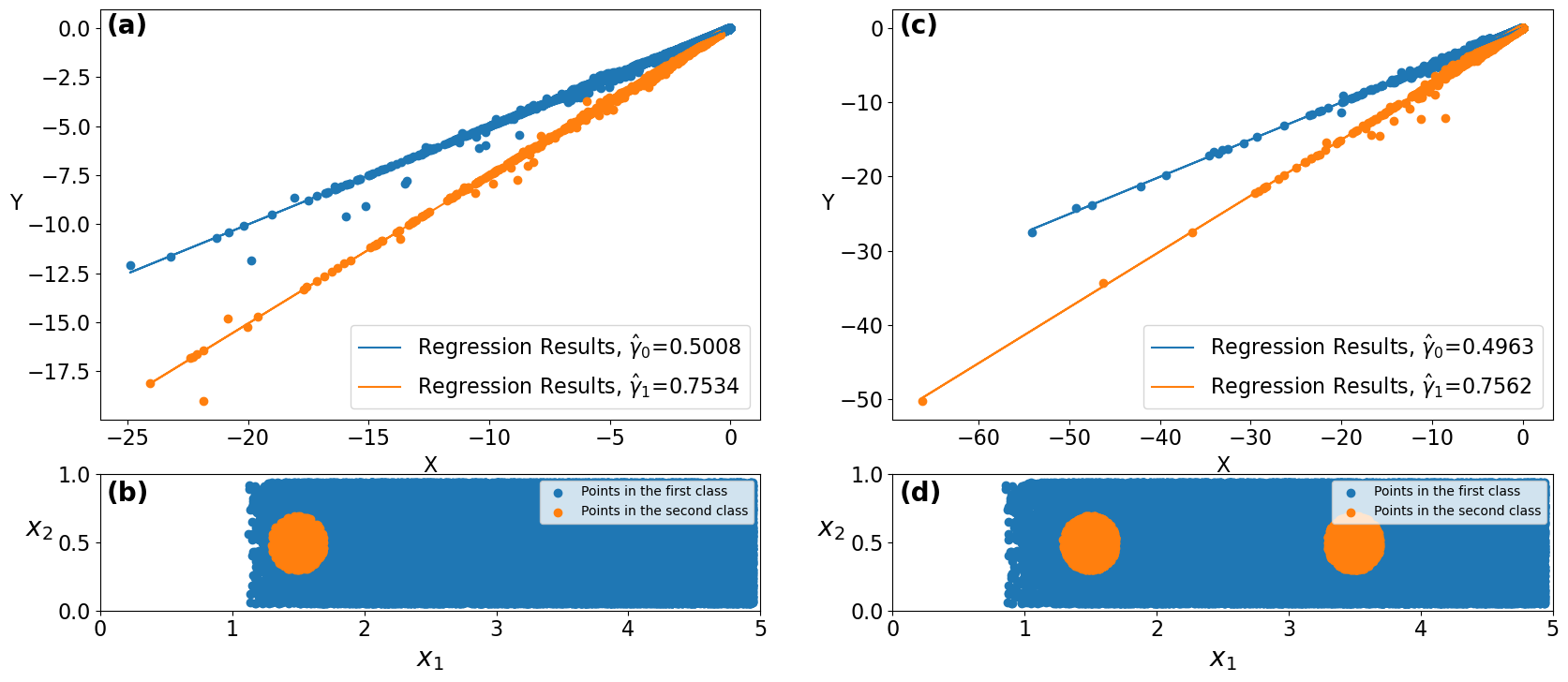}
				\end{centering}
				\caption{
					\NYY{
						Propagation in inhomogeneous media with a disk
						inclusion after the
						classification using the $k$-means algorithm.
						The solid line is the result of
						linear regression for each class.
						(b) The positions of the data point
						corresponding to each class. 
						(c) The sampled data in $(X,Y)$ for the crack
						propagation in inhomogeneous media with two disk
						inclusions after the
						classification using the $k$-means algorithm.
						The solid line is the result of
						linear regression for each class.
						(d) The positions of the data point
						corresponding to each class. 
					}
					\NYYY{
						The applied strain $A$ in \eqref{BC_strain} is chosen as $A=1.25$ for
						(a,b) and $A=1$ for (c,d).
					}		 			 
				}  \label{1disk_2diskscase1_difclasses_AT2}
			\end{figure}
			
			Our method successfully estimate the toughness of
			inclusions at different positions.
			In the cases (VI) and (VII) in
			Fig.\ref{Fig_Allpaths_AT2}, the crack
			bypasses the inclusions (VI), or it bypass one inclusions
			and is stuck at the second inclusion.
			Due to the bypassing, the number of data points
			corresponding to the inclusions is
			small.
			Still, we may estimate the toughness as $\hat{\gamma}_0=0.4994$ and
			$\hat{\gamma}_1=0.7493$ for the case (VI), and $\hat{\gamma}_0=0.4988$ and
			$\hat{\gamma}_0=0.7475$ for the case (VII).
			\NYb{
				We may also estimate the position of the inclusions as
				shown in Fig.\ref{2hardercase2_difclasses} in
				Supplementary Materials.
			}
		}		
		
		\NYb{
			We should stress that we use the data before the crack reaches the
			second inclusion in the case of two disk inclusions (see
			Fig.~\ref{fig_crack_tip_at2}(v)-(vii) with the time interval $[100\Delta t, 200 \Delta t)$.
			Therefore, we may {\it predict} the position and
			toughness of the second inclusion.
		}
		
		\subsection{inverse problems of the $\mathrm{AT_1}$ model}
		
		\NY{
			Next we consider the phase-field model
			\eqref{Phase_field_model_Bourdin}, which is based on the
			surface energy $\mathrm{AT_1}$.
			Naively, we may apply the same algorithm as in
			section \ref{section_Inverse_problem_Takaishi_Kimura}
			and change the corresponding definitions in step 4, the \textbf{Linear Regression}.
			However, as we will see, this naive method does not work
			for the $\mathrm{AT_1}$ model.
			This is because in this model, the plus operator
			$(\cdot)_+$ is effective almost everywhere in space
			except near the crack tip.
			Therefore, we need additional pre-processes to estimate
			the toughness in this model.
			
			We assume that the values of $\alpha_2$, $\epsilon$, and $\mu$ are known.  
			Based on the second equation of system
			\eqref{Phase_field_model_Bourdin},  if we omit the plus
			operator $(\cdot)_+$, formally we obtain
			\begin{equation}\label{def_estimate_Bourdin}
				\alpha_2  \frac{\partial \widetilde{Z}}{\partial t}
				=\gamma \left(\frac{3}{4}\epsilon \Delta \widetilde{Z}
				- \frac{3}{8\epsilon}\right) +\mu \lvert
				\nabla \widetilde{U}\rvert^2(1-\widetilde{Z}).
			\end{equation}
			This leads us to define 
			\begin{equation*}\label{def_Y_Bourdin}
				Y:= \alpha_2  \frac{\partial \widetilde{Z}}{\partial t} - \mu \lvert \nabla \widetilde{U}\rvert^2 (1-\widetilde{Z}) 
			\end{equation*}
			and 
			\begin{equation}\label{def_X_Bourdin}
				X :=  \frac{3}{4}\epsilon \Delta \widetilde{Z}-\frac{3}{8\epsilon}.
			\end{equation}
		}
		If we find the linear relation such that $Y=
		\hat{\gamma} X$, then $\hat{\gamma}$
		is the estimator of $\gamma$, that is the fracture
		toughness of the phase-field model.

		
		\NY{\subsubsection{homogeneous case}\label{section_Bourdin_homogeneous}
			
			We first perform the estimation of the homogeneous case
			by using the same algorithm as in the Section~\ref{section_Inverse_problem_Takaishi_Kimura}.
			The sampled data is shown in Fig.~\ref{TotalSpace_homo_AT1} in the Supplementary Materials.
			The deviation is due to the plus operator $(\cdot)_+$.
			In fact, when the right-hand side of
			\eqref{def_estimate_Bourdin} is negative, the crack
			field $\widetilde{Z}$ does not satisfy \eqref{def_estimate_Bourdin}
			but $\partial \widetilde{Z}/\partial t=0$.
			This occurs particularly when $\widetilde{Z}$ and
			$\widetilde{U}$ are uniform in space, corresponding to the
			points away from the crack tip.
			The difference from the $\mathrm{AT_2}$ originates from
			the second term in \eqref{def_X} and
			\eqref{def_X_Bourdin}.
			When $\widetilde{Z} \approx 0$, $X \approx 0$ in
			\eqref{def_X} and accordingly the right-hand side of
			\eqref{def_estimate_TakaishiKimura} is approximately 0.
			On the other hand, in the $\mathrm{AT_1}$ when
			$\widetilde{Z} \approx 0$, $X < 0$ and the right-hand side
			of \eqref{def_estimate_Bourdin} becomes negative if
			$\widetilde{U}$ is approximately uniform in space.
			In this case, $X$ and $Y$ do not exhibit a linear
			relationship due to the plus operator.
			
			

			
			We, therefore, propose an additional pre-process of the
			data.
			Inspired by the form of the solution in
			Fig.~\ref{Solution_Cut_AT2}, we propose to sample the
			points at the region where the values of $z$ are closer to $1$.			
			Given a time interval $[n_0\Delta t, n_1 \Delta t)$, at
			each time step, we consider the region where the crack
			has not arrived, similar to the estimation for the
			{$\mathrm{AT_2}$ model}.
			We denote the tip position at $n\Delta
			t$ to be $(x_{1,n}, x_{2,n})$ and define $\Omega_n =
			[x_{1,n},x_{1,n}+\Delta x]\times [x_{2,n},x_{2,n}+\Delta y]$, that is a small disk containing four
			points around the crack tip. 
			As we use the uniform square mesh, the sample points
			correspond to one-fourth of a disk whose radius $r$
			satisfies $(\sqrt{5}/2)\Delta x<r<(3/2)\Delta
			x$.

			\begin{figure}[htbp] 
				\begin{center}
					\includegraphics[scale=0.3]{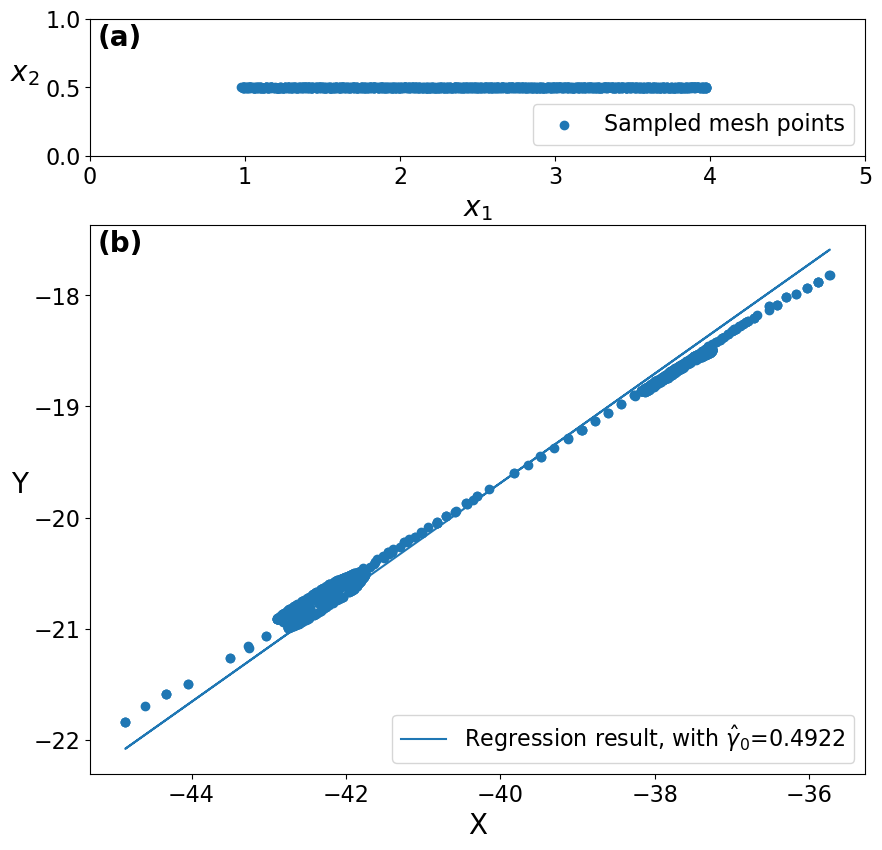}
				\end{center}
				\caption{
					\NYY{
						(a) The positions of the data point obtained
						by the $\mathrm{AT_1}$ model after pre-processing. 			 
						(b) The sampled data in $(X,Y)$ for the crack
						propagation in homogeneous media using the $\mathrm{AT_1}$ model.
						The solid line is the result of
						linear regression.
					}
					\NYYY{
						The applied strain $A$ in \eqref{BC_strain} is chosen as $A=1.25$.
					}		 			 
				} \label{Total_XY_homo_tip_Bourdin_theoretical}
			\end{figure}
					%

			We consider the time interval $[n_0\Delta t, n_1 \Delta t)$.
			The $j^{th}$ sample point will be denoted by $p_j:=((x_{1}^{(j)},x_{2}^{(j)}),t^{(j)})$ and the total sample set is denoted by $P =\{p_j\}_{j=1,2,...,N_{data}}$ with $N_{data} = (n_1-n_0)\cdot 2 \cdot 2$. 
			The sampled data points of the time interval $[80\Delta t, 380 \Delta t)$ for the homogeneous media are
			shown in \NYY{Fig.~\ref{Total_XY_homo_tip_Bourdin_theoretical}(a).}
			We note that the crack path is a straight line in the
			homogeneous material case (see Fig.\ref{Fig_Allpaths_AT2}).
			In contrast with the estimation of $\mathrm{AT_2}$
			model, the number of data points is small and limited
			to the region near the crack path.
			Therefore, the estimation of the toughness and positions
			of inclusions are much harder in $\mathrm{AT_1}$ model.
			
			We compute the $X$ and $Y$ on the chosen data set.
			The result is shown in \NYY{Fig.~\ref{Total_XY_homo_tip_Bourdin_theoretical}(b)}.
			By applying the linear regression, we can estimate
			\NYY{$\hat{\gamma}=0.4922$}, which is comparable with the
			ground truth $\gamma=0.5$.
			Despite the \NYb{small number of} sample points, our method successfully
			estimate the toughness of the homogeneous media in the $\mathrm{AT_1}$ model.

			
			\subsubsection{Inhomogeneous case}\label{section_Bourdin_Inhomogeneous}
			
			For the inhomogeneous materials in the
			$\mathrm{AT_1}$, we apply the $k$-means as in
			Sec.~\ref{section_Inverse_problem_Bourdin} to
			divide the data into two classes corresponding
			to two toughness $\gamma_0$ and $\gamma_1$.
			We also have to use the pre-processing the data
			points as in Sec.~\ref{section_Bourdin_homogeneous}.
			We consider the time interval to be $[80\Delta t, 380 \Delta
			t)$ for the cases (II) and (III), $[80\Delta t, 300 \Delta
			t)$ for the case (IV) and $[80\Delta t, 480 \Delta
			t)$ for the cases (V)-(VII).
			
			\NYb{
				Figure~\ref{Stripes_points_Bourdin_greendisk}(a) shows
				the data after pre-processing for the data of
				crack propagation in the inhomogeneous media
				with the stripe inclusions.} We may see the data can be decomposed into two classes.
			In fact, by applying $k$-means, we find two classes
			corresponding to two toughness, $\gamma_0$ in the
			homogeneous media and $\gamma_1$ in the inclusions.
			The regression results show \NYY{$\hat{\gamma}_0=0.5037$ and
				$\hat{\gamma}_1=0.751$}.
			They are reasonably closed to the ground truth
			$\gamma_0=0.5$ and $\gamma_1=0.75$.
			
			The data points in each class are located at $\gamma_0$
			or $\gamma_1$ as shown in \NYY{Fig.~\ref{Stripes_points_Bourdin_greendisk}(b)}.
			We are able to discover the $x_1$ position
			of the inhomogeneous part.
			However, in contrast with the estimation shown in
			Sec.~\ref{section_Inverse_problem_Takaishi_Kimura}, we
			cannot identify all the position of
			the inhomogeneity.
			This is because the effect of the plus operator is
			stronger in the $\mathrm{AT_1}$ model than the
			$\mathrm{AT_2}$ model.
			We cannot use the data points at which the plus operator
			works.
			Therefore, we lose many data points in the
			$\mathrm{AT_1}$ model.
			Still, we may find the position of tougher regions using
			our estimation method.				
			
			\begin{figure}[htbp] 
				\begin{center}
					\includegraphics[scale=0.3]{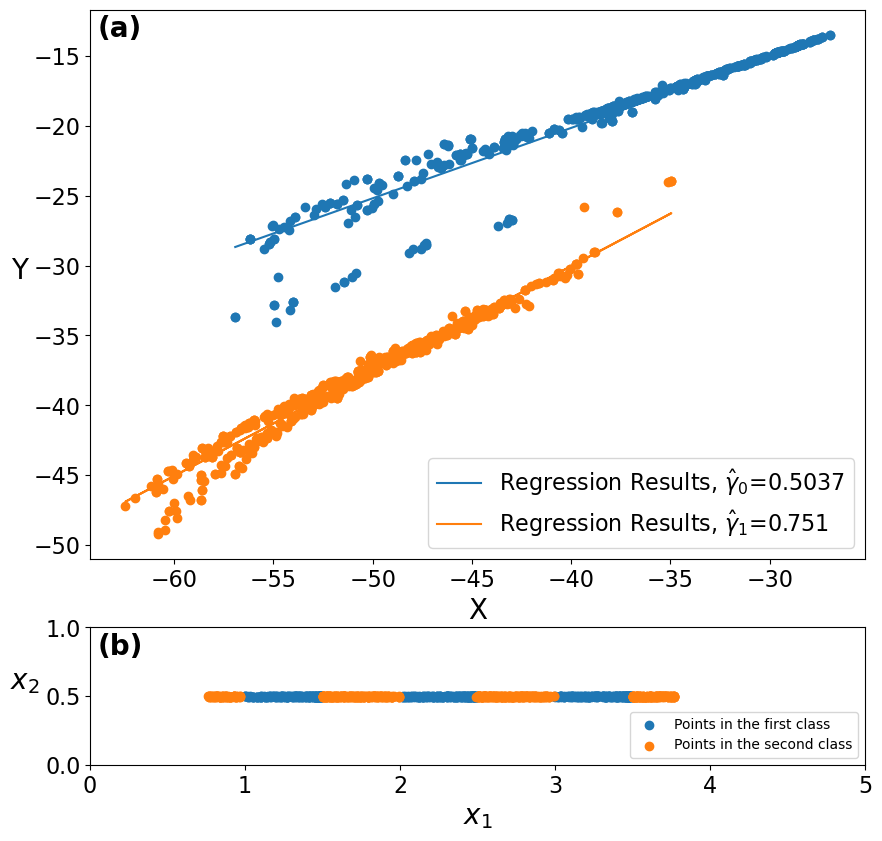}
				\end{center}
				\caption{
					\NYY{
						(a) The sampled data in $(X,Y)$ obtained by the
						$\mathrm{AT_1}$ model for the crack
						propagation in inhomogeneous media with stripe
						inclusions after the
						classification using the $k$-means algorithm.
						The solid line is the result of
						linear regression for each class.
						(b) The positions of the data point
						corresponding to each class. 
					}
					\NYYY{
						The applied strain $A$ in \eqref{BC_strain} is chosen as $A=1.25$.
					}		 			 
				}  \label{Stripes_points_Bourdin_greendisk}
			\end{figure}
			
			\NYY{Figure~\ref{Disk_points_Bourdin}(a)} shows the results of
			estimation for one-disk inclusion.		
			The regression results show \NYY{$\hat{\gamma}_1=0.4934$ and  $\hat{\gamma}_2=0.7471$}. 
			\NYY{Figure~\ref{Disk_points_Bourdin}(b)} shows the points of
			$\gamma_0$ and $\gamma_1$ in the \NYb{$x_1-x_2$ space}.
			We can estimate the position of $x_1$ coordinates of the inhomogenous part, that is the disk inclusion.
			We recall the the inclusion is a disk with barycenter $(1.5,0.5)$ and radius 0.2.
			%
			\NYb{
				The estimation also works for the one harder disk inclusion in
				the $\mathrm{AT_1}$ model as shown in
				Fig.\ref{Backto_spacecoordinates_disk1_Bourdin}
				in Supplementary Materials.
			}
			
			\begin{figure}[htbp] 
				\begin{centering}
						\includegraphics[scale=0.2]{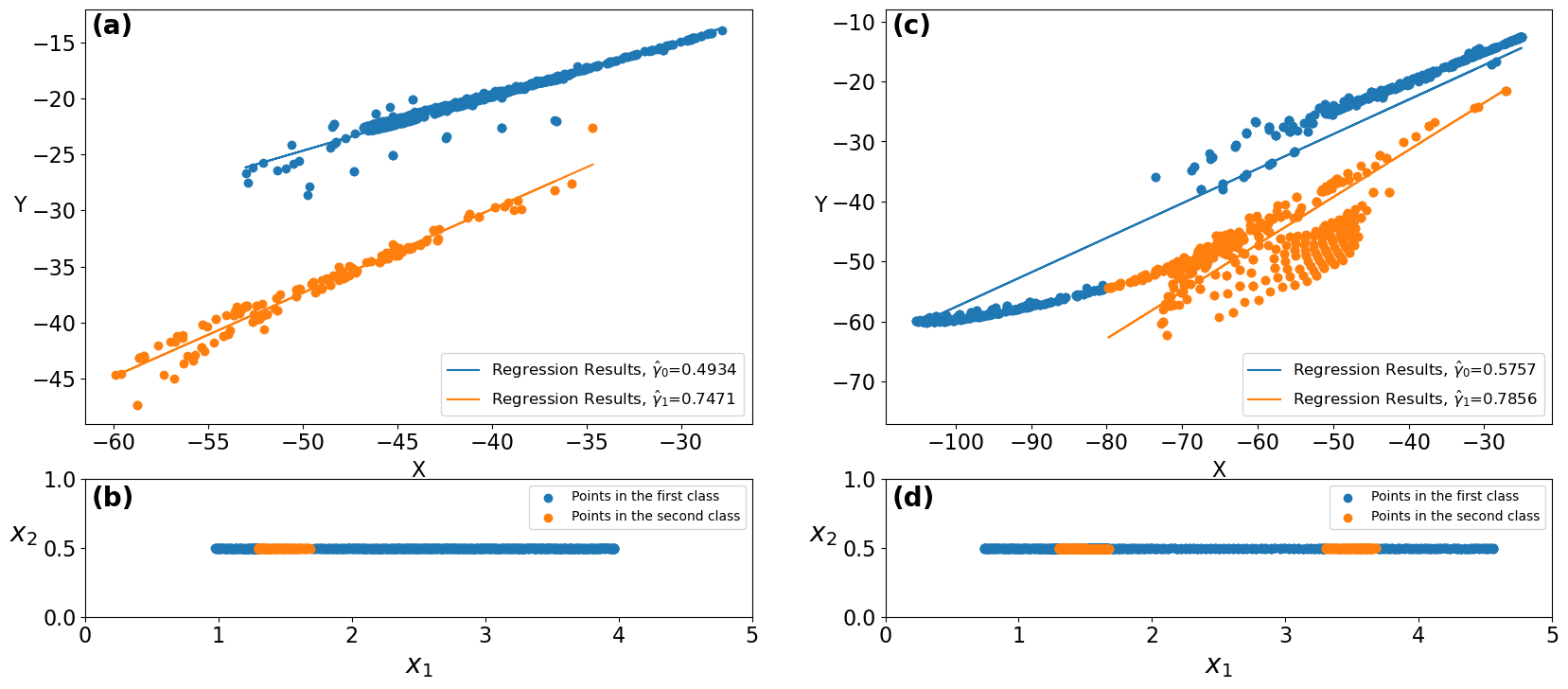}
					\end{centering}
					\caption{
						\NYY{
							(a) The sampled data in $(X,Y)$ obtained by the
							$\mathrm{AT_1}$ model for the crack
							propagation in inhomogeneous media with a disk
							inclusion after the
							classification using the $k$-means algorithm.
							The solid line is the result of
							linear regression for each class.
							(b) The positions of the data point
							corresponding to each class. 
							(c) The sampled data in $(X,Y)$ obtained by the
							$\mathrm{AT_1}$ model  for the crack
							propagation in inhomogeneous media with two disk
							inclusions.
							The data is classified using the $k$-means algorithm.
							The solid line is the result of
							linear regression for each class.
							(d) The positions of the data point
							corresponding to each class. 
						}
						\NYYY{
							The applied strain $A$ in \eqref{BC_strain} is chosen as $A=1.25$ for
							(a) and $A=1$ for \NYb{(c)}.
						}		 			 
					}  \label{Disk_points_Bourdin}
				\end{figure}
				
				\NYb{
					Similar to
					Sec.\ref{section_Inverse_problem_Bourdin},
					we consider the two-disk inclusions.
					As in Fig.~\ref{Disk_points_Bourdin}(c,d), the estimation is less accurate
					compared with the $\mathrm{AT_2}$ model.
					Still, we estimate the toughness as $\hat{\gamma}_1=0.5757$
					and $\hat{\gamma}_1=0.7856$, which are
					reasonably close to the ground truth
					$\gamma_0=0.5$ and $\gamma_1=0.75$.
					Even when the two inclusions are placed at
					the different $x_2$ position
					(Fig.\ref{Fracture_Toughness_6cases}(VI,VII)),
					we may estimate the toughness of the
					inclusions and their positions (see Fig.~\ref{2hardercase2_difclasses_AT1} in Supplementary Materials).
					In the $\mathrm{AT_1}$ model, because we
					use the data near the crack tip, the
					information is limited to the spatial
					points near the crack path.
					Therefore, we cannot estimate the whole
					shape of the inclusions.
					We may also need data close to the
					inclusions.
					This is in contrast with the $\mathrm{AT_2}$ model.
				}
			}

			\section{Discussion}\label{section_discussion}
			
			\NY{
				In summary, we study numerical simulations of the crack
				propagation using the phase-filed models, and propose an
				algorithm to estimate the fracture toughness as well as
				the position of the inhomogeneity.
				The method is based on pre-processing of the data and
				linear regression.
				One important step is to sample the grid points at each
				time in the region where the crack has not arrived to
				avoid the effect of the plus operator.
				We successfully demonstrate to estimate the fracture
				toughness of the inclusion and its position.
				The method works both for strip and disk inclusions, and
				both for one- and two-disk inclusions.
			}
			
			\NYYY{			
				In this work, the two different types of energy functional,
				the $\mathrm{AT_1}$ and $\mathrm{AT_2}$ models,
				are studied.
				Both models exhibit the increase of the
				$J$-integral corresponding to the inhomogeneous
				toughness.
				For the inverse problems, our method works for
				both models, but the $\mathrm{AT_1}$ model
				requires more pre-processing than the
				$\mathrm{AT_2}$ model.
				Moreover, the estimation of the position of
				inclusions is limited for the $\mathrm{AT_1}$
				model, but not for the $\mathrm{AT_2}$ model.
				This is because the plus operator is effective
				in a broader
				region in the $\mathrm{AT_1}$ model.
				In this respect, the $\mathrm{AT_2}$ model gives
				us more information about the inclusions in the
				inverse problems.
				Nevertheless, the $\mathrm{AT_1}$ model has an
				advantage that the crack field $z({\bf x},t)$
				has a sharper distribution near the crack, and
				strictly $z=0$ away from the crack.
				\NYb{
					This difference may affect the precision
					of the computation of the $J$-integral
					during the crack propagation.
				}
			}
			
			\NYYY{
				We should note that our method of the estimation of
				inhomogeneous toughness is based on parameter
				estimation, and not on model discovery, which is a
				recent active research
				field \cite{Rudy_Brunton_Proctor_Kutz_2017}.
				In the model discovery, a model described by partial
				differential equations is not {\it a priori} known;
				ideally, elastic equation and phase-field model are
				estimated from the data.
				\NYb{
					We assume that we know the phase-field
					model and the linear elastic equation.
					There are several difficulties in the model discovery in
					our system.
				}
				First, we consider slow quasi-static crack propagation,
				and therefore, $\partial_t z \approx 0$.
				The model discovery based on the regression method
				relies on minimization of the cost function describing
				the difference between the left-hand side and right-hand
				side of the equation with sparse regularization \cite{Rudy_Brunton_Proctor_Kutz_2017}. 
				When the left-hand side ($\partial_t z$) vanishes, we cannot
				uniquely estimate all the coefficients in the right-hand
				side.
				Second difficulty is the large parameter contrast in the
				model of 
				\eqref{Phase_field_model_TakaishiKimura} and
				\eqref{Phase_field_model_Bourdin} due to $\epsilon$.
				The regression with sparse regularization tries to
				remove a term whose coefficient is small.
				However, when the model contains a term with larger
				coefficient $\mathcal{O}(1/\epsilon)$ and one with
				smaller coefficient $\mathcal{O}(\epsilon)$, it is
				difficult to keep the $\mathcal{O}(\epsilon)$ term, but
				remove other unnecessary terms.
			}
			
				\GY{
				For the crack propagation in a homogeneous media, after
				using the same pre-processing that we use, we may apply the SINDy
				method, proposed in \cite{Rudy_Brunton_Proctor_Kutz_2017}.
				If we consider the dictionary terms from
				the product between $z$ and $u$ and their derivatives: such as $\{u_x, u_y, u_x^2, u_y^2, z, z_{x},z_{y}, z_{xx},z_{yy}, zu_x, zu_y, zu_x^2, zu_y^2, zz_{x},zz_{y}, zz_{xx}, zz_{yy}\}$,
				we can estimate
				unnecessary terms particularly when the tolerance is
				large, (corresponding to weaker sparse
				regularization).				
				This is due to the large parameter contrast.
				The first issue mentioned above is not too serious.
				If we include only the ground-truth terms in the
				dictionary, our estimation result is
				\begin{equation}\label{Estimated_SINDy_Giventerms}
					z_t = 0.52447 \Delta z -1057.1z -42.221 \lvert \nabla u \rvert^2
					+42.350 z\lvert \nabla u \rvert^2,
				\end{equation}		
				under the ground-truth parameter $\gamma=0.5$, $\mu=1$, and $\epsilon=0.02$.
			}
			The estimation does not look working, but the ratio
			between the first and third terms, and between the second and fourth term in \eqref{Estimated_SINDy_Giventerms} are $\hat{\epsilon} \hat{\gamma}/\hat{\mu} =0.012$ and $ \hat{\gamma}/(\hat{\epsilon} \hat{\mu}) =24.96$, respectively.
			Therefore, the estimation works up to constant
			multiplication.
			
			\NYYY{
				The application of SINDy for the crack propagation in
				inhomogeneous media is more involved.
				Using the same assumption made in our study, namely, if we assume only
				two values of toughness, we may use the group
				sparsity \cite{Rudy:2018}.
				Nevertheless, the difficulties discussed above 
				apply also in the inhomogeneous toughness, and
				therefore, the estimation of the model is still
				impractical.
				
				In the original version of SINDy, the estimation is not
				robust against noise, and accurate computations of
				spatial and time derivatives are
				necessary \cite{Schaeffer2017}.
				In our method, we use the Chebyshev polynomial of
				degree {7}.
				We have tried different interpolation schemes to check
				how the accuracy of spatial and time derivatives affects
				the estimation.
			}
			\NYb{
				For example, in the stripe case of the $\mathrm{AT_2}$ model, by
				applying the second-order finite
				difference method to compute the
				derivatives of $\widetilde{U}$ and
				$\widetilde{Z}$,  we obtain the
				estimation $\hat{\gamma}_0=0.5015$ and
				$\hat{\gamma}_1=0.75$. The result is
				shown in Fig.~\ref{Stripes_AT2_FiniteDifference} in
				Supplementary Materials.			
				Even in such a poor evaluation of derivatives, our
				estimation gives the values close to the
				ground truth, suggesting that our method does
				not require very accurate computations of spatial and
				time derivatives.
				We speculate this is because we assume we {\it a priori}
				know the model, and therefore, the search space during
				the optimization is limited.
			}
			

			\NYY{
				In our method, we can compute the strain field $\{U_K^n\}$ in \eqref{Phase_field_model_TakaishiKimura} and \eqref{Phase_field_model_Bourdin} once the
				crack field $\{Z_K^n\}$ is known.
				This means that we may estimate inhomogeneous toughness only from the
				information of a crack path.
				This aspect of the estimation is important because in practical
				situations, it is difficult to measure a strain field.
				The strain field should be obtained or estimated from the crack field.
				However, our method is limited because we assume the linear elastic
				model and also the elastic constant is known.
				If these assumptions are not available, we need to estimate the model
				for the elastic equation, and even for the linear elastic model, we
				need to estimate the elastic constant.
				In this case, the estimation only from data of a crack field is more
				involved.
				A possible extension of our method is to treat the strain field as a
				hidden variable, and estimate the inhomogeneous toughness after
				marginalizing the hidden variable \cite{Bishop_2006}.
			}


			%

			\section*{Appendix A. $k$-means algorithm}\label{section_kmeans} 
			In the inhomogeneous case, we used an extension of the classic $k$-means classification to our problem. The algorithm is as follows:
			\begin{enumerate}
				\item We compute the maximum slope and minimum slope of the sample points and denote the values by $a_{max}$ and $a_{min}$ respectively. Then, we choose $\{\gamma_m\}_{m=1,2,...k}$ that is $k$ independent random values which follow the uniform distribution $\mathcal{U}(\gamma_{min}, \gamma_{max})$; 
				\item For the sample points $P=\{p_j\}$, we define $$\operatorname{dist}_m(p_j)=\left[Y(p_j) -\gamma_m X(p_j)\right]^2$$ and assigned $p_j$ to class $\mathrm{arg\,min}_{m=1,2,...k}\operatorname{dist}_m(p_j)$. We assign all the sample points into $k$ classes.
				\item In each class, we perform linear regression to update the value of $\{\gamma_m\}_{m=1,2,...k}$;
				\item We go back to step 2 and reassign the sample points to $k$ classes and then to step 3 to update the values of $\{\gamma_m\}_{m=1,2,...k}$. We repeat the algorithm till that $\{\gamma_m\}_{m=1,2,...k}$ do not change anymore.
			\end{enumerate}

			\section*{Acknowledgement}
			\NYb{
				The research work of Y.G. was supported
				by the research funding of MathAM-OIL, AIST c/o AIMR, Tohoku
				University. The authors acknowledge  the support by JSPS KAKENHI grant numbers Number JP19K14605 to
				Y.G. and JP20K03874 to N.Y..
			}
			
			
			\bibliographystyle{elsarticle-num}
			
			
			\newpage
				\section*{Supplement Material}
			\subsection{Additional data of the estimation of homogeneous material in $\mathrm{AT_1}$}
			\begin{figure}[htbp] 
				\begin{center}
					\includegraphics[scale=0.3]{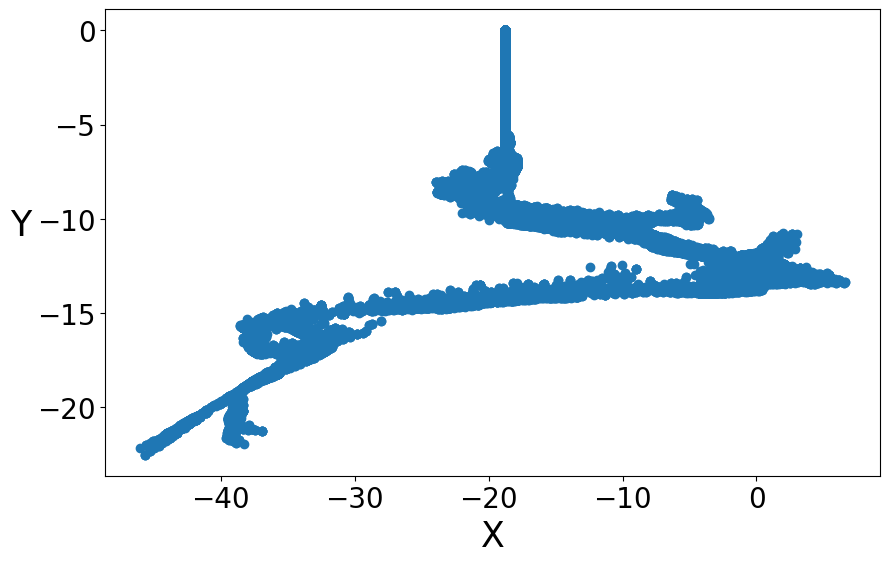}
				\end{center}
				\caption{
					The sampled data in $(X, Y)$ for the crack propagation in homogeneous media using the $\mathrm{AT}_1$ model, using the same algorithm for the $\mathrm{AT}_2$ model.
				}  \label{TotalSpace_homo_AT1}
			\end{figure}
			
			\subsection{Additional data of the estimation of inhomogeneous toughness}
			
			\begin{figure}[htbp] 
				\begin{center}
					\includegraphics[scale=0.2]{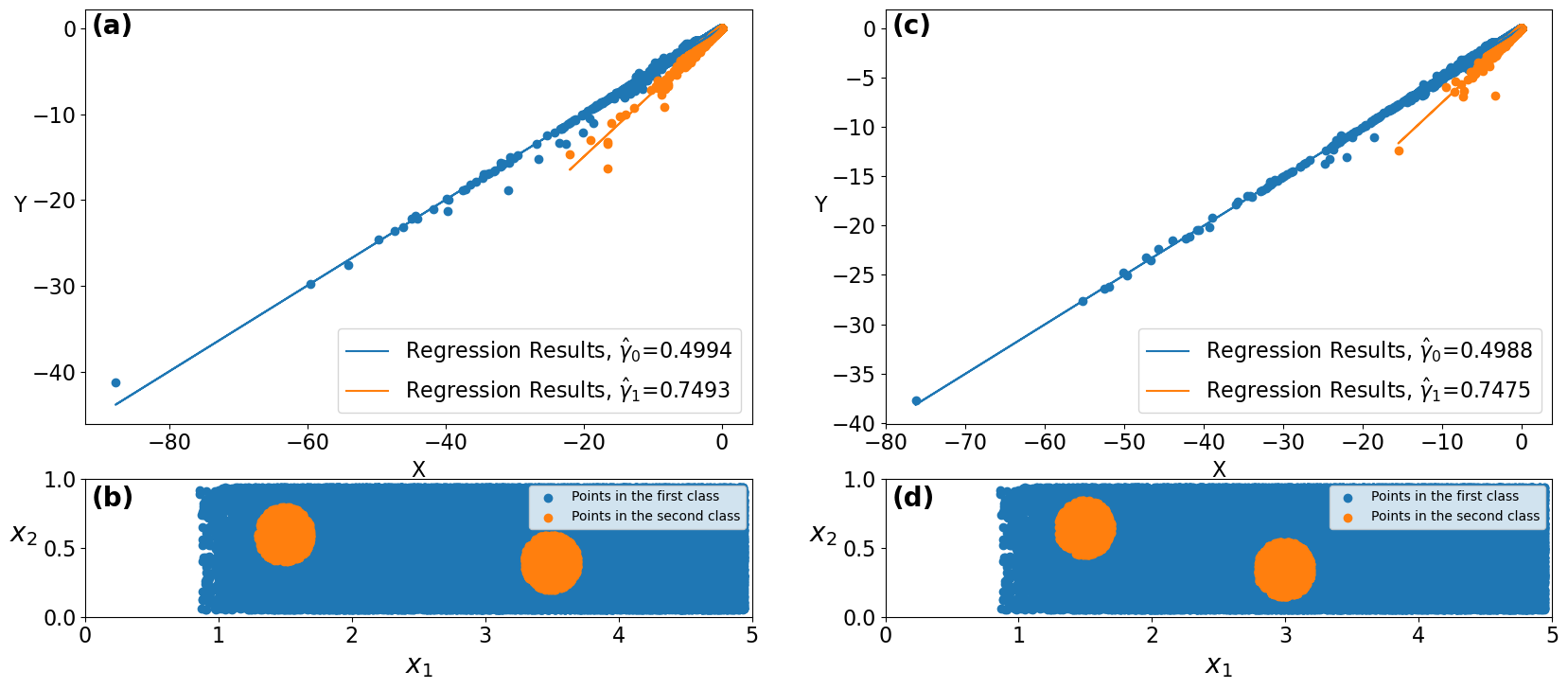}
				\end{center}
				\caption{
					\NYY{
						(a) The sampled data in $(X,Y)$ obtained by the
						$\mathrm{AT_2}$ model for the crack
						propagation in inhomogeneous media with two disk
						inclusions.
						The data is classified by using the $k$-means algorithm.
						The solid line is the result of
						linear regression for each class.
						(b) The positions of the data point
						corresponding to each class. 
						(c) The sampled data in $(X,Y)$ obtained by the
						$\mathrm{AT_2}$ model for the crack
						propagation in inhomogeneous media with two disk
						inclusions.
						The data is classified by using the $k$-means algorithm.
						The solid line is the result of
						linear regression for each class.
						(d) The positions of the data point
						corresponding to each class. 
					}
					\NYYY{
						The applied strain $A$ in \eqref{BC_strain} is chosen as $A=1$.
					}		 		 
				}  \label{2hardercase2_difclasses}
			\end{figure}

			\begin{figure}[htbp] 
				\begin{center}
					\includegraphics[scale=0.3]{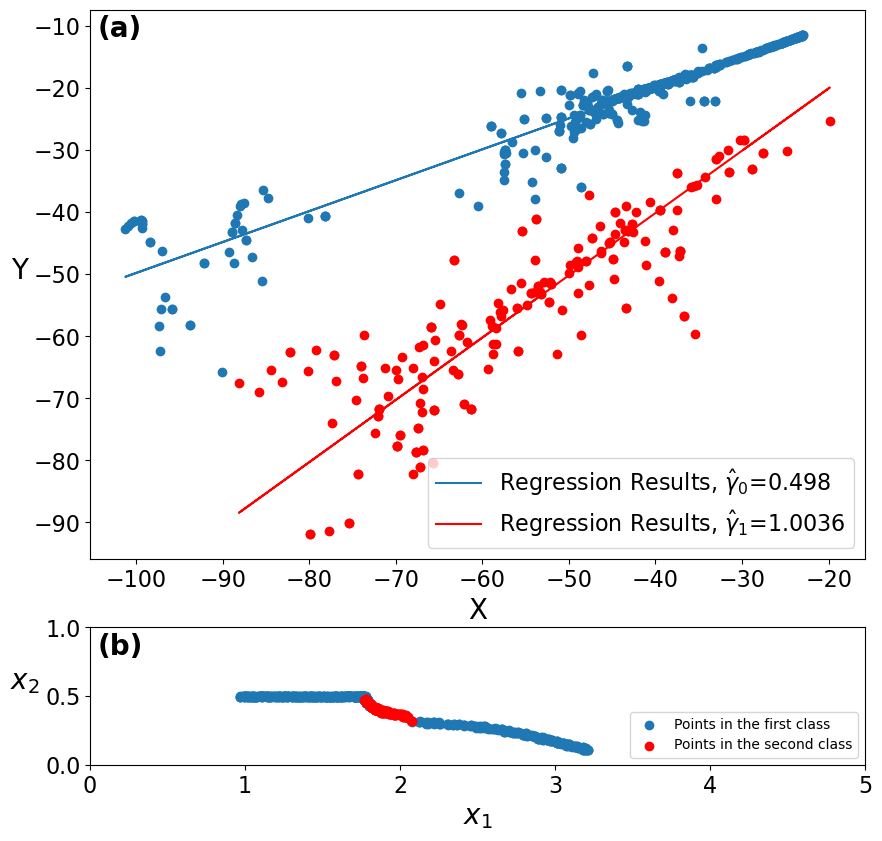}
				\end{center}
				\caption{
					\NYY{
						(a) The sampled data in $(X,Y)$ obtained by the
						$\mathrm{AT_1}$ model in the time interval $[80\delta t, 300 \delta t)$ for the crack
						propagation in inhomogeneous media with a
						harder disk
						inclusion.
						The data is classified by using the $k$-means algorithm.
						The solid line is the result of
						linear regression for each class.
						(b) The positions of the data point
						corresponding to each class. 
					}
					\NYYY{
						The applied strain $A$ in \eqref{BC_strain} is chosen as $A=1.25$.
					}		 		 
				}  \label{Backto_spacecoordinates_disk1_Bourdin}
			\end{figure}

			\begin{figure}[htbp] 
				\begin{center}
					\includegraphics[scale=0.2]{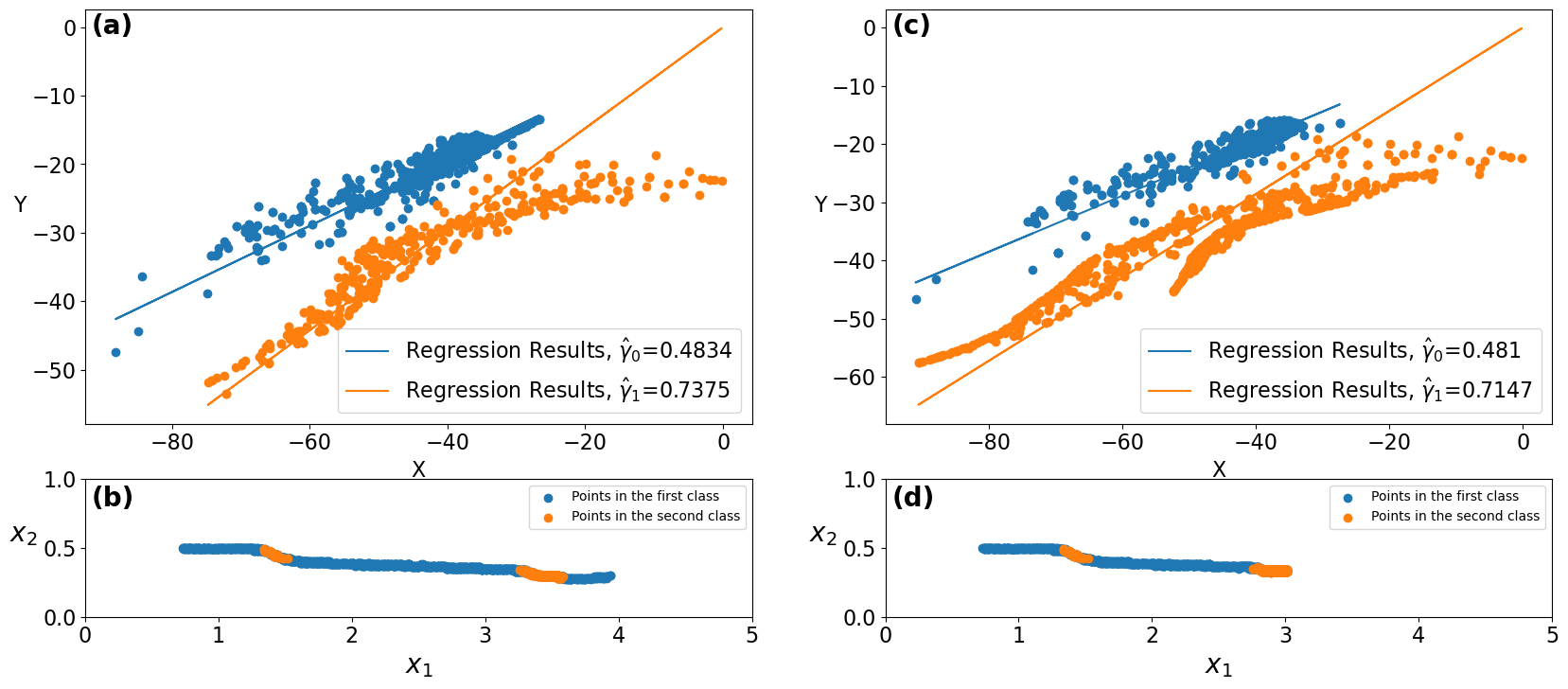}
				\end{center}
				\caption{
					\NYY{
						(a) The sampled data in $(X,Y)$ obtained by the
						$\mathrm{AT_1}$ model in the time interval $[80\delta t, 480 \delta t)$  for the crack
						propagation in inhomogeneous media with two disk inclusions.
						The data is classified by using the $k$-means algorithm.
						The solid line is the result of
						linear regression for each class.
						(b) The positions of the data point
						corresponding to each class. 
						(c) The sampled data in $(X,Y)$ obtained by the
						$\mathrm{AT_1}$ model in the time interval for the crack
						propagation in inhomogeneous media with two disk inclusions.
						The data is classified by using the $k$-means algorithm.
						The solid line is the result of
						linear regression for each class.
						(d) The positions of the data point
						corresponding to each class. 
					}
					\NYYY{
						The applied strain $A$ in \eqref{BC_strain} is chosen as $A=1$.
					}		 		 
				}  \label{2hardercase2_difclasses_AT1}
			\end{figure}

			\begin{figure}[htbp] 
				\begin{center}
					\includegraphics[scale=0.3]{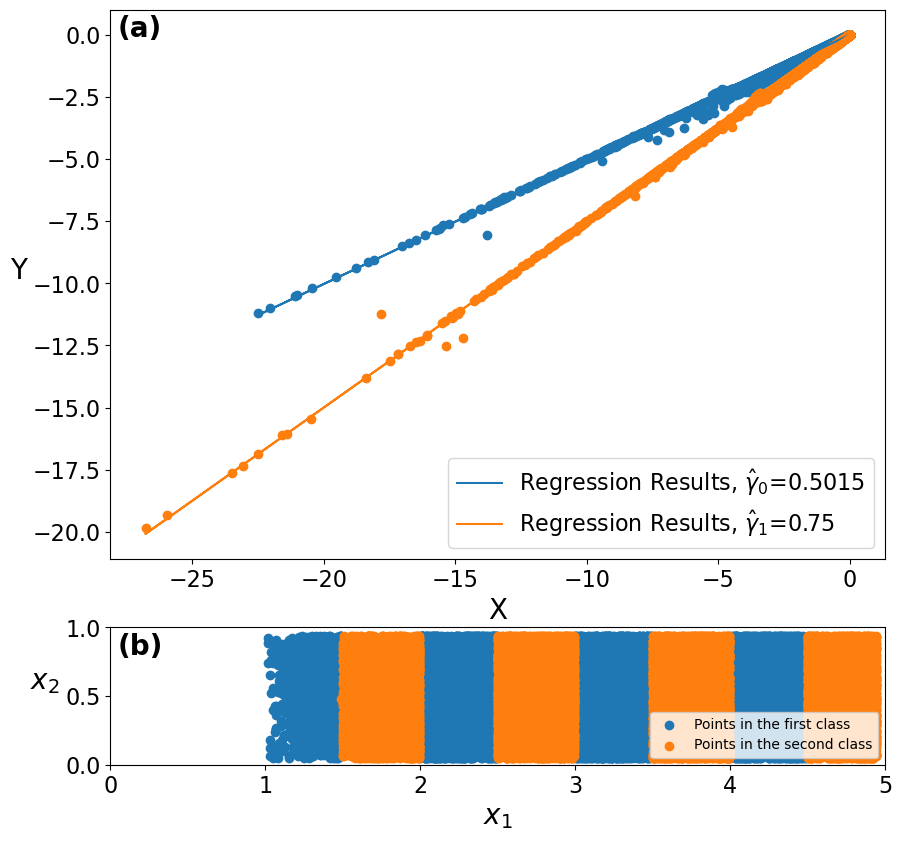}
				\end{center}
				\caption{
					{
						(a) The sampled data in $(X,Y)$ obtained by the
						$\mathrm{AT_2}$ model for the crack
						propagation in inhomogeneous media with stripe inclusion. The derivatives are computed by the finite difference method. 
						The data is classified by using the $k$-means algorithm.
						The solid line is the result of
						linear regression for each class.
						(b) The positions of the data point
						corresponding to each class. 
					}
					\NYYY{
						The applied strain $A$ in \eqref{BC_strain} is chosen as $A=1.25$.
					}		 		 
				}  \label{Stripes_AT2_FiniteDifference}
			\end{figure}
			
		\end{document}